\documentclass[12pt]{amsart}

\usepackage{float,wrapfig}

\numberwithin{equation}{section}

\usepackage{dsfont,amsfonts,amsmath,pgf}

\newcommand{\eps}{\varepsilon}

\newcommand{\RR}{\mathds{R}}

\newcommand{\NN}{\mathds{N}}

\usepackage{fancybox}
      \newtheorem{theorem}{Theorem}[section]
       \newtheorem{proposition}[theorem]{Proposition}
       \newtheorem{corollary}[theorem]{Corollary}
       \newtheorem{lemma}[theorem]{Lemma}
  
\theoremstyle{remark}
       \newtheorem{remark}{Remark}[section]
         
\theoremstyle{definition}
\newtheorem{definition}{Definition}[section]

\def\E{{\mathbb E}}
\def\V{{\rm Var}}
\def\Cov{{\rm Cov}}
\def\cov{\Cov}

\newcommand{\calF}{{\mathcal F}}

\newcommand{\wo}{\mathbb{E}}                                                    
\newcommand{\wwo}[2]{\mathbb{E}\left[\left.{#1}\right|{#2}\right]}              
\DeclareMathOperator{\Var}{Var}                                                 
\newcommand{\wVar}[2]{\Var\left[\left.{#1}\right|{#2}\right]}                   
\DeclareMathOperator{\sign}{sign}                                               
\newcommand{\indyk}[1]{1\!\!\!\! 1_{#1}}                                        

\newcommand{\TT}{\mathcal{T}}

\def\Z{\mathbf Z}

 \newcommand{\ceta}{\eta}
 \newcommand{\ctheta}{\theta}

 \usepackage{framed}
\colorlet{shadecolor}{gray!30}
\newenvironment{ana}
  {\begin{leftbar}
  \begin{shaded} }
{  \end{shaded}\end{leftbar}}

 \newcommand{\arxiv}[1]{\begin{ana}
  #1\end{ana}}

\date{Created:  October 27, 2009. \\ Printed: \today \ file: \jobname.tex}

\usepackage[top=.8 in,bottom=.8 in, left=.8 in, right=.8 in]{geometry}

\author[W.~Bryc]{
W{\l}odek  Bryc
}
\address{
Department of Mathematical Sciences,
University of Cincinnati,
PO Box 210025,
Cincinnati, OH 45221--0025, USA}
\email{Wlodzimierz.Bryc@UC.edu}

\author[W.~Matysiak]{Wojciech Matysiak}
\address{Wydzia{\l} Matematyki i Nauk Informacyjnych\\
Politechnika Warszawska\\
Pl. Politechniki 1\\
00-661 Warszawa, Poland}
\email{matysiak@mini.pw.edu.pl}

\title[Wilson's $6-j$ laws and Markov processes]{Wilson's $6-j$ laws and stitched Markov processes}
\begin{document}
 \maketitle
 \begin{abstract}
We show how to insert time into the parameters of the Wilson's $6-j$ laws to construct
discrete Markov chains with these laws. By a quadratic transformation we convert them into Markov  processes with linear regressions and quadratic conditional variances. Further conversion into the "standard form" gives  "quadratic harnesses"  with  "classical" value  of parameter $\gamma=1\pm 2\sqrt{\sigma\tau}$.
For   $\gamma=1+2\sqrt{\sigma\tau}$,   a random-parameter-representation of the original Markov
 chain allows us to stitch together two copies of the process, extending %
 time domain of the quadratic harness %
   from $(0,1)$ to $(0,\infty)$.
  \end{abstract}
\arxiv{This  is an expanded version with additional details that are omitted from the version intended for publication.}

\section{Introduction}
The work on this paper started with an attempt to fit Markov processes with linear regressions and quadratic conditional variances into Wilson's $6-j$-laws from \cite{Wilson:1980}. This required choosing appropriate time-parameterization of the laws so that we get a Markov chain, and the appropriate (quadratic) transformation  of this chain so that  conditional and absolute moments are given by simple enough formulas.

Generically, processes with linear regressions and quadratic conditional variances can be further transformed ("standardized") so that they are described by five parameters, see  \cite[Theorem 2.2]{Bryc-Matysiak-Wesolowski-04}. We  expected Wilson's $6-j$ laws to lead to the "classical" quadratic harnesses with the parameters tied by equality $\gamma=1-2\sqrt{\sigma\tau}$. But, to our surprise, depending on the range of parameters we also got  quadratic harnesses with $\gamma=1+2\sqrt{\sigma\tau}$. In the latter case, the initial construction gave only a quadratic harness with time  $(0,1)$. However,  the underlying Markov chain is a mixture of simpler Markov chains. We used this mixture representation to extend the quadratic harnesses to $(0,\infty)$  by stitching together  two conditionally-independent Markov chains with shared randomization. The stitching approach was suggested by the construction of the "bi-Pascal" process with $\gamma=1+2\sqrt{\sigma\tau}$ in \cite{Maja:2009}; %
our argument is modeled on \cite{Bryc-Wesolowski-10}.

The paper is organized as follows. In Section \ref{Sect:QHFNV} we use Wilson $6-j$ laws to construct quadratic harnesses on $(0,1)$ or on $(0,\infty)$, depending on the range of parameter $C$. These are Case 1 and Case 2 of Theorem \ref{Cor-parameters}.
In Section \ref{Sect:EQHC2} we represent Markov chain from Section \ref{Sect:QHFNV} as a mixture of "simpler" Markov chains. We also confirm that each of these Markov chains transforms into a quadratic harness with $\gamma=1$ and $\sigma=0$ (which is our justification for the adjective "simpler" in the previous sentence.) In Section \ref{Sect:stitch} we stitch together a pair of such quadratic harnesses into the quadratic harness on $(0,\infty)$, thus extending the process from Case 1   of Theorem \ref{Cor-parameters} to the maximal %
time domain.

The expanded version of this paper with additional technical or computational details is posted on the arXiv.
\subsection{Quadratic harnesses}\label{S: QH} %
In \cite{Bryc-Matysiak-Wesolowski-04} the authors  consider
 square-integrable  stochastic processes on $(0,\infty)$
such that
for all $t,s> 0$,
\begin{equation}\label{EQ:cov}
\E(Z_t)=0,\: \E(Z_s Z_t)=\min\{s,t\},
\end{equation}
and for $s<t<u$, $\E({Z_t}|{\mathcal{F}_{ s, u}})$ is a linear function of $Z_s,Z_u$, and
$\V [Z_t|\mathcal{F}_{s,u }]$ is a quadratic  function of $Z_s,Z_u$.
Here, $\mathcal{F}_{ s, u}$ is the two-sided
$\sigma$-field generated by $\{Z_r: r\in (0,s]\cup[u,\infty)\}$.
Then  \eqref{EQ:cov} implies that
\begin{equation}
\label{EQ:LR} \E({Z_t}|{\mathcal{F}_{ s, u}})=\frac{u-t}{u-s}
Z_s+\frac{t-s}{u-s} Z_u
\end{equation}
 for all $s<t<u$, which is sometimes referred to as a harness condition, see e.~g. \cite{Mansuy-Yor04}.

While there are numerous examples of harnesses, the  assumption of quadratic conditional variance is more restrictive. For  example, all integrable L\'evy processes are harnesses, but as determined
 by Weso\l owski \cite{Wesolowski93}, only a few of them are also quadratic harnesses.
Under certain  technical assumptions,  \cite[Theorem 2.2]{Bryc-Matysiak-Wesolowski-04}
 asserts that quadratic variance has the following form: there exist numerical constants
 $\eta,\theta,\sigma,\tau,\gamma\in\RR$ such that
for all $s<t<u$,
\begin{multline}\label{EQ:q-Var}
\V [Z_t|\mathcal{F}_{s,u }]\\
= \frac{(u-t)(t-s)}{u(1+\sigma s)+\tau-\gamma s}\left( 1+\eta \frac{u Z_s-s Z_u}{u-s} +\theta\frac{Z_u-Z_s}{u-s}\right. \\
 \left.
+ \sigma
\frac{(u Z_s-s Z_u)^2}{(u-s)^2}+\tau\frac{(Z_u-Z_s)^2}{(u-s)^2}
-(1-\gamma)\frac{(Z_u-Z_s)(u Z_s-s Z_u)}{(u-s)^2} \right).
\end{multline}

\begin{definition}\label{DEF-harn}
We will say that a square-integrable stochastic process $(Z_t)_{t\in T}$ is a quadratic harness on $T$ with parameters $(\eta,\theta,\sigma,\tau,\gamma)$,
if it satisfies \eqref{EQ:cov}, \eqref{EQ:LR} and \eqref{EQ:q-Var} on an open interval $T$ which may be all of or a proper subset of $(0,\infty)$.  We also assume that the one-sided conditional moments are as follows:
 for $0< s\leq t\leq u$ in $T$,
\begin{equation}\label{LR+}
\wwo{Z_t}{\calF_{\geq u}}=\frac{t}{u}Z_u,
\end{equation}
\begin{equation}\label{LR-}
\wwo{Z_t}{\calF_{\leq s}}=Z_s,
\end{equation}
\begin{eqnarray} \label{CV-}
\wVar{Z_t}{\calF_{\leq s}}&=&\frac{t-s}{1+\sigma s}\left(\sigma Z_s^2+
\eta Z_s+1\right),\\ \label{RV-}
\wVar{Z_t}{\calF_{\geq u}}&=&\frac{t(u-t)}{u+\tau}\left(\tau \frac{Z_u^2}{u^2}+\theta \frac{Z_u}{u}+1\right).
\end{eqnarray}
\end{definition}
We remark that on infinite intervals, formulas (\ref{LR+}-\ref{RV-}) follow from the other assumptions, see \cite[(2.7), (2.8), (2.27), and (2.28)]{Bryc-Matysiak-Wesolowski-04}.

We expect that quadratic harnesses on finite intervals are determined uniquely by the parameters.  This has been confirmed under some technical assumptions when the parameters satisfy additional constraints, of which the main constraint seem to have been that $-1\leq \gamma\leq 1-2\sqrt{\sigma\tau}$.

 It is known, see \cite{Bryc-Matysiak-Wesolowski-04},  that for quadratic harnesses on $(0,\infty)$, parameters $\sigma,\tau$ are non-negative, and that
 $ \gamma \leq 1+2\sqrt{\sigma\tau}$. Quadratic harnesses with $\gamma=1- 2\sqrt{\sigma\tau}$ were called "classical" in \cite{Bryc-Matysiak-Wesolowski-04}. Quadratic harnesses with $\gamma=1+2\sqrt{\sigma\tau}$ could also have been called "classical", but
 there had been no examples of such processes until the bi-Pascal process was constructed in \cite{Maja:2009}.  The bi-Pascal process does not have higher moments so large part of general theory developed in \cite{Bryc-Matysiak-Wesolowski-04} does not apply.
 Our interest here is in providing additional examples  of quadratic harnesses with $\gamma=1\pm 2\sqrt{\sigma\tau}$.

\section{Quadratic harnesses with finite number of values}\label{Sect:QHFNV}
A family of quadratic harnesses  $(Z_t)$ with two values appears in
\cite[Section 3.2]{Bryc-Matysiak-Wesolowski-04b}. These processes have parameter $\gamma=-1$ and their trajectories follow two quadratic curves. Since $1,Z_t,Z_t^2$ are linearly dependent, the parameters in \eqref{EQ:q-Var} are not determined uniquely. In fact, one can show that for this family of processes the admissible parameters in \eqref{EQ:q-Var}  can take any real values $\sigma,\tau\in\RR$ (positive or negative) such that $\sigma\tau\ne 1$  and  any $\eta,\theta\in\RR$ such that
\begin{equation*}
(\eta\tau+\theta )(\eta+\theta\sigma) +\left(1-\sigma\tau\right)^2>0.
\end{equation*}

  Quadratic harnesses with finite number of values, including two values,  appear also in  \cite[Section 4.2]{Bryc-Wesolowski-08}. The processes constructed there have parameter $\gamma<1-\sqrt{\sigma\tau}$.

In this section we construct (non-homogeneous) Markov processes which take
a finite number of values and  we show how to transform them into quadratic harnesses with $\gamma=1\pm2\sqrt{\sigma\tau}$.
 These processes are different from the previous ones even in the case of two-values; this can be seen from analyzing the curves they follow, see Figure \ref{F1}. Somewhat surprisingly,  processes corresponding to $\gamma=1\pm\sqrt{\sigma\tau}$ are described by the same formulas for transition probabilities, differing only in the range of one of the parameters that enter the formulas.

Our construction is based on Wilson's \cite{Wilson:1980}  $6-j$ laws. As in \cite[(3.5)]{Wilson:1980}
we fix integer $N\geq 1$  and assume that
\begin{equation}\label{A1}
a>-1/2,\; b\in(-a,a+1), \mbox{ and either $c>a+N$ or
$c<-a-N+1$}.
\end{equation}
(The choice of the range for $c$ will later affect the properties of the quadratic harness.)

 For  $N\in\NN$ and $k=0,1,\dots,N$, define
 \begin{equation}\label{p_k}
p_{k,N}(a,b,c)=C_N(a,b,c)
\frac{(2 a)_k (a+1)_k (a+b)_k (a+c)_k (-N)_k}{k! (a)_k (a-b+1)_k (a-c+1)_k (2 a+N+1)_k}.
\end{equation}
where the normalizing constant is
\begin{equation}\label{K}
C_N(a,b,c)=\frac{(a-b+1)_N (a-c+1)_N}{(2 a+1)_N (-b-c+1)_N}
\end{equation}
and
$(a)_k=\Gamma(a+k)/\Gamma(a)=a(a+1)\dots (a+k-1)$ is the Pochhammer symbol.

When the parameters are such that numbers $p_{k,N}(a,b,c)$ are well defined, then the sum over all $k=0,\dots,N$ is one; this is \cite[formula (3.4)]{Wilson:1980} applied to $m=n=0$.  So under assumption \eqref{A1}, from \cite[(3.4)]{Wilson:1980} one reads out that
\begin{equation}\label{nu}
\nu(dx)=\sum_{k=0}^N p_{k,N}(a,b,c) \delta_k(dx)
\end{equation}
  is a probability measure on $\{0,1,\dots,N\}$.

The following algebraic formula will be used several times. %
\begin{lemma}\label{L.2.1}
For $k=0,1,\dots,N$, $j=0,1,\dots,k$,
\begin{equation}\label{M-P}
\frac{p_{j,k}(a,b,a+k+2 \delta) p_{k,N}(a+\delta,b+\delta,c)}{p_{j,N}(a,b,c+\delta) }=
p_{k-j,N-j}(a+j+\delta,-a-j+\delta,c)\;.
\end{equation}

\end{lemma}
\begin{proof}
Multiplying both sides of \eqref{M-P} by $p_{j,N}(a,b,c+\delta)$, expanding them by the use of \eqref{p_k} and \eqref{K}, canceling out common terms and grouping the remaining ones, we observe that \eqref{M-P} would follow if we verify that $LHS=RHS$ with
\begin{multline*}
LHS=
\frac{1}{(2a+1)_k(2a+k+1)_j}
\frac{(2a+2\delta)_k(2a+2\delta+k)_j}{(2a+2\delta+1)_N(2a+2\delta+N+1)_k}
(a+c+\delta)_k \\
\frac{(a+\delta+1)_k}{(a+\delta)_k}
\frac{(-2\delta-k+1)_k}{(-2\delta-k+1)_j}
\frac{(a-c+\delta+1)_N}{(a-c+\delta+1)_k}
\frac{(a+b+2\delta)_k}{(-a-b-2\delta-k+1)_k}
\frac{(-k)_j(-N)_k}{k!},
\end{multline*}
and
\begin{multline*}
RHS=
\frac{(2a+2j+1)_{N-j}}{(2a+1)_N(2a+N+1)_j(2a+2j+1)_{k-j}}\\
\frac{(2a+2\delta+2j)_{k-j}}{(2a+2\delta+2j+1)_{N-j}(2a+2\delta+j+N+1)_{k-j}}
(a+c+\delta)_j(a+c+\delta+j)_{k-j}\\
\frac{(a+\delta+j+1)_{k-j}}{(a+\delta+j)_{k-j}}
(2\delta)_{k-j}
\frac{(a-c+\delta+j+1)_{N-j}}{(a-c+\delta+j+1)_{k-j}}\\
\frac{(a-c-\delta+1)_N}{(a-c-\delta+1)_j(a-c-\delta+j+1)_{N-j}}\frac{(j-N)_{k-j}(-N)_j}{(k-j)!}.
\end{multline*}

To perform the verification, we will use the following simplification rules
\begin{eqnarray}
(\alpha)_N(\alpha+N)_M &=& (\alpha)_{M+N},\label{jeden}\\
\frac{(\alpha+L+1)_{M-L}}{(\alpha+L)_{M-L}} &=& \frac{\alpha+M}{\alpha+L}\label{dwa},\\
\frac{(\alpha)_L}{(\alpha+1)_{M+L}} &=& \frac{\alpha}{(\alpha+L)_{M+1}}\label{trzy}.
\end{eqnarray}

From \eqref{jeden} it follows that
\begin{multline*}
\frac{(2a+2j+1)_{N-j}}{(2a+1)_N(2a+N+1)_j(2a+2j+1)_{k-j}}=\frac{(2a+k+j+1)_{N-k}}{(2a+1)_{N+j}}=\\
\frac{(2a+k+j+1)_{N-k}}{(2a+1)_{k+j}(2a+1+k+j)_{N-k}}=
\frac{1}{(2a+1)_k(2a+k+1)_j}.
\end{multline*}
Similarly, \eqref{jeden} and \eqref{trzy} give
\begin{multline*}
\frac{(2a+2\delta+2j)_{k-j}}{(2a+2\delta+2j+1)_{N-j}(2a+2\delta+j+N+1)_{k-j}}\\=
\frac{(2a+2\delta+2j)_{k-j}}{(2a+2\delta+2j+1)_{N+k-2j}}=
\frac{2a+2\delta+2j}{(2a+2\delta+k+j)_{N+1-j}},
\end{multline*}
while applied to the analogous expression in $LHS$ they give
\begin{multline*}
\frac{(2a+2\delta)_k(2a+2\delta+k)_j}{(2a+2\delta+1)_N (2a+2\delta+N+1)_k}=
\frac{(2a+2\delta)_{k+j}}{(2a+2\delta+1)_{k+N}}=
\frac{2a+2\delta}{(2a+2\delta+k+j)_{N+1-j}}.
\end{multline*}
From \eqref{jeden} we get
\[
(a+c+\delta)_j (a+c+\delta+j)_{k-j}=(a+c+\delta)_k,
\]
\[
\frac{(a-c+\delta+j+1)_{N-j}}{(a-c+\delta+j+1)_{k-j}}=(a-c+\delta+k+1)_{N-k}=\frac{(a-c+\delta+1)_N}{(a-c+\delta+1)_{k}},
\]
and (this expression appears only in $RHS$)
\[
\frac{(a-c-\delta)_N}{(a-c-\delta+1)_j(a-c-\delta+j+1)_{N-j}}=1.
\]
From \eqref{dwa} it follows that
\[
\frac{(a+\delta+j+1)_{k-j}}{(a+\delta+j)_{k-j}}=\frac{a+\delta+k}{a+\delta+j}\ \textrm{ and }\
\frac{(a+\delta+1)_k}{(a+\delta)_k}=\frac{a+\delta+k}{a+\delta}.
\]
For the remaining expressions from $LHS$ we have
\[
\frac{(-2\delta-k+1)_k}{(-2\delta-k+1)_j}=(-1)^{k+j}(2\delta)_{k-j}\ \mbox{ and }\ \frac{(a+b+2\delta)_k}{(-a-b-2\delta-k+1)_k}=(-1)^k.
\]

Now, the above simplifications show that equality $LHS=RHS$ is equivalent to
\begin{multline*}
\frac{(-k)_j(-N)_k}{k!}(-1)^{k+j}(2\delta)_{k-j}\frac{2a+2\delta}{(2a+2\delta+k+j)_{N+1-j}}\frac{a+\delta+k}{a+\delta}(-1)^k\\=
\frac{(j-N)_{k-j}(-N)_j}{(k-j)!}(2\delta)_{k-j} \frac{2a+2\delta+2j}{(2a+2\delta+k+j)_{N+1-j}}\frac{a+\delta+k}{a+\delta+j},
\end{multline*}
which is easily seen to be true. Thus \eqref{M-P} is proved.

\end{proof}
In particular, by taking the sum over $j$ in \eqref{M-P} we have
\begin{equation}\label{M-P-2}
p_{k,N}(a+\delta,b+\delta,c)=\sum_{j=0}^k p_{j,N}(a,b,c+\delta) p_{k-j,N-j}(a+j+\delta,-a-j+\delta,c).
\end{equation}

Next we compute the moments of an auxiliary random variable associated with  probability law \eqref{nu}.

\begin{proposition}\label{P:M+V} %
 Suppose parameters $a,b,c,N$ satisfy \eqref{A1}.
For $k=0,\dots,N$, consider a random variable $Y$ such that
\begin{equation}\label{Y}
\Pr(Y=k(2a+k))=p_{k,N}(a,b,c).
\end{equation}
Then
\begin{equation}\label{EX-WilsonD}
\E(Y)=\frac{(a+b) (a+c) N}{b+c-N}
\end{equation}
 and
\begin{equation}\label{VX-WilsonD}
\Var(Y)=
\frac{N (a-b+N) (a+b) (a-c+N) (a+c) (b+c)}{(b+c-N)^2 (N-b-c-1)}.\end{equation}
\end{proposition}
\begin{proof}The proof is elementary for $N=0,1$, as the law of $Y$ is  $\delta_0$ for $N=0$ and
 $$
\frac{(a-b+1) (a-c+1)}{(2 a+1) (-b-c+1)} \delta_{0}+\frac{  (a+b) (a+c)}{(2 a+1)  (b+c-1)}\delta_{2a+1}
 $$
 for $N=1$.

For $N>1$ and  $k\geq 1$, we have
$$k(2a+k)p_{k,N}(a,b,c)=\frac{(a+b) (a+c) N}{b+c-N}p_{k-1,N-1}(a+1,b,c) \;,$$
which gives \eqref{EX-WilsonD}. For $k>2$, iterating the algebraic identity we get
 \begin{multline*}
 k(k-1)(2a+k)(2a+k+1)p_{k,N}(a,b,c)\\
 =\frac{(a+b) (a+b+1) (a+c) (a+c+1) (N-1) N}{(b+c-N) (b+c-N+1)}p_{k-2,N-2}(a+2,b,c).
 \end{multline*}
 Noting that $k (k - 1) (2 a + k) (2 a + k + 1)=k^2 (2 a + k)^2  - (1 + 2 a) k (2 a + k)$, we  evaluate
 $$\E(Y^2)=\frac{(a+b) (a+b+1) (a+c) (a+c+1) (N-1) N}{(b+c-N) (b+c-N+1)}+(1+2a)\E(Y)$$ and \eqref{VX-WilsonD} follows  by an elementary calculation.
\end{proof}

\subsection{Markov chain}\label{Sect:Markov}
Now we introduce %
a continuous time (non-homogeneous) Markov chain on the finite state space $\{0,1,\dots,N\}$   with  parameters   $A$, $B$, $C$. We assume that $N\in\NN$,  $A>-1/2$, $B\in(-A,A+1)$. For the third parameter,   we will assume that either
\begin{itemize}
\item[\bf Case 1:] $C<-A-N+1$,
\end{itemize}
or
\begin{itemize}
\item[\bf Case 2:] $C>A+N$.
\end{itemize}
These two cases will appear in several statements below.

The Markov process will be defined for  $t\in \TT$, where
\begin{equation*}
\TT=\begin{cases}
\left(-(A+B),\infty\right) &\mbox{  in Case 1,}\\
\left(-(A+B),C-A-N\right) &\mbox{ in Case 2.}
\end{cases}
\end{equation*}
We remark that $A+B>0$ and that in Case 2 the interval $\TT$ is non-empty, as $C-A-N>0$.

 We  also remark that the process in Case 1 is well defined on another interval $(-\infty,C-A-N)$. %
 This second "component" of the process will be used to extend the quadratic harness from $(0,1)$ to $(0,\infty)$. (In fact, one should think that in Case 1, the process starts at $-A-B$ at state $0$, continues through $\infty=-\infty$ and ends at state $N$ at time $C-A-N$.)

For $s<t$ in $\TT$, define %
matrix  $P_{s,t}=[p_{s,t}(k,n)]_{0\leq k,n\leq N}$ with entries
\begin{equation}\label{xiK-tr}
p_{s,t}(k,n)=p_{n-k,N-k}\left(A+\frac{t}{2}+k,-A-s+\frac{t}{2}-k,C-\frac{t}{2}\right)
\end{equation}
if $0\leq k\leq n \leq N$, and let $p_{s,t}(k,n)=0$ for all other values $k,n\in\{0,\dots,N\}$.

The following shows that matrices $P_{s,t}$  are transition probabilities
of a Markov chain.
\begin{proposition}[Chapman-Kolmogorov equations]
For $s<t<u$ in $\TT$, %
$p_{s,t}(k,n)\geq 0$, and
for $j\leq n\leq N$, 
\begin{equation}\label{CK-2}
p_{s,u}(j,n)=\sum_{k=j}^n p_{s,t}(j,k)p_{t,u}(k,n).
\end{equation}
Furthermore, %
\begin{equation}\label{CK-3}
\sum_{k=j}^N p_{s,t}(j,k) =1.
\end{equation}

\end{proposition}
\begin{proof}
To verify \eqref{CK-2}, we apply \eqref{M-P-2} with parameters
$$
a=A+j+t/2, \; b=-A-j-s+t/2 ,\; c=C-u/2 ,\; \delta= (u-t)/2.
$$
Formula \eqref{CK-3} is the already mentioned generic identity for the weights.

To verify that  for %
$0\leq k\leq n\leq N$
we have $p_{s,t}(k,n)\geq 0$ we verify assumption \eqref{A1}.
Here $a=A+k+t/2$, $b=-A-k-s+t/2$, $c=C-t/2$. We get
$$a=A+k+t/2>(A-B)/2+k>-1/2+k>0$$
as $t>-(A+B)$, $B<A+1$ and $k\ge0$,
$$b=-A-k-s+t/2=-a+t-s>-a$$
as $t> s$, and
$$a+1-b=2A+1+2k+s>A-B+2k+1>2k>0$$
as $s>-(A+B)$, $B<A+1$ and $k\geq0$. If $C<-A-N+1$ then
$$c=C-t/2<-A-N+1-t/2 %
=-a-(N-k)+1.$$ If $C>A+N$ then
$$c-a=C-A-k-t>N-k,$$
since $t<C-A-N$.

\end{proof}

Let $(\xi_t)_{t\in \TT}$ be the Markov chain constructed above, i.e.  %
$$\Pr(\xi_t=n|\xi_s=k)=p_{s,t}(k,n),\, 0\leq k\leq n\leq N.
$$
In particular, the univariate laws of the Markov chain are
$$\Pr( \xi_t=j)=p_{-A-B,t}(0,j)=p_{j,N}(A+t/2,B+ t/2,C-t/2).$$

For $t\in \TT$, consider the process
\begin{equation}\label{Y-ini}
Y_t=(2A+t+\xi_{t})\xi_t+A(A+t) = (A+t+\xi_{t})(A+\xi_{t}).
\end{equation}

\begin{lemma}\label{L.2.3}
$(Y_t)$ is a Markov process with mean
\begin{equation}\label{m}
\E(Y_t)=\frac{A(B+C)(A+N)+N B C }{B+C-N}+t\frac{A (B+C)+N C}{B+C-N}
\end{equation}
and variance
\begin{equation}\label{v}
\Var(Y_t)=-\frac{N(A+C) (B+C) (A-B+N) }{(B+C-N)^2 (B+C-N+1)}(A+B+t) (A-C+N+t).
\end{equation}
\end{lemma}

\begin{proof}
The Markov property follows from the fact that the lines $\ell_k(t)= (A+t+k)(A+k)$ for $k=0,1,\dots, N$ do not intersect over $t\in \TT$. Indeed, $\ell_k$ and $\ell_j$ intersect at $t=-(2A+j+k)<-2A-1<-(A+B)$. Therefore, the law $\sum p_t(k)\delta_k$ converges to the degenerate law $\delta_{0}$ as $t\to -(A+B)$ from the right.

The formulas for the mean and the variance are now recalculated from Proposition \ref{P:M+V}, noting that $Y_t$ is in distribution $Y+A(A+t)$ with $Y$ given by \eqref{Y}.
\end{proof}
Of course, $(Y_t)_{t\in\TT}$  naturally extends to the left endpoint by $\lim_{t\to -(A+B)}Y_t = -AB $ in mean square and, for a separable version, almost surely.

\pgfdeclareimage[width=3in]{T}{traj}
\begin{figure}[H]
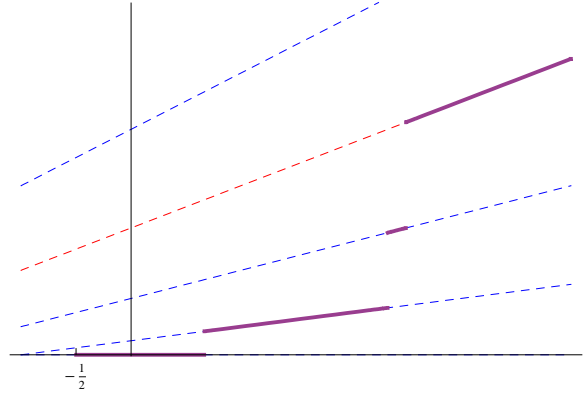

\pgfuseimage{T}
\caption{Trajectories of process $(Y_t)_{t\in\TT}$ lie on the family of   lines.
The process starts at line $\ell_0 =(A+t)A$, at $t=-A-B$ and jumps up
until it stops on a random  line  %
 $\ell_\Theta=(A+t+\Theta)(A+\Theta)$
that it follows to $\infty$. %
(Here, $A=0, B=1/2, N=4, \Theta=3$; for more on $\Theta$, see \eqref{K-cond}.)
\label{F1}
}
\end{figure}

From Proposition \ref{P:M+V} we read out that for $s<t$
the conditional moments are
\begin{equation}\label{cond1}
\E(Y_t|Y_s)=Y_s\frac{A-C+N+t}{A-C+N+s} -\frac{C(A+N)(t-s)}{A-C+N+s},
\end{equation}
\begin{multline}\label{vcond1}
\Var(Y_t|Y_s)\\=(A-C+N+t)(t-s)\frac{\left(Y_s+C(s-C)\right)\left(Y_s-(A+N)(s+A+N)\right)}{(A-C+N+s)^2(A-C+N+s-1)}.
\end{multline}
Indeed, $p_{s,t}(k,n)=p_{n-\xi_s,N-\xi_s}(A+t/2+\xi_s,-A-s+t/2-\xi_s,C-t/2)$.
Recall
\begin{equation}\label{G:cov-g}
\Cov(X,Y)=\E(\Cov(X,Y|U))+\Cov(\E(X|U),\E(Y|U)).
\end{equation}
Using \eqref{G:cov-g} with $X=U=Y_s$, $Y=Y_t$, from \eqref{v} and \eqref{cond1} we compute 
\begin{multline}\label{cov}
\Cov(Y_s,Y_t)=\frac{A-C+N+t}{A-C+N+s} \Var(Y_s)
\\=-\frac{N(A+C) (B+C) (A-B+N) (A+B+s) (A-C+N+t)}{(B+C-N)^2 (B+C-N+1)}.
\end{multline}

We now compute the two-sided conditional distribution
\begin{equation*}
\Pr(\xi_t=j|\xi_s=k, \xi_u=m)=\frac{\Pr(\xi_t=j|\xi_s=k)\Pr(\xi_u=m|\xi_t=j)}{\Pr(\xi_u=m|\xi_s=k)}
=\frac{p_{s,t}(k,j)p_{t,u}(j,m)}{p_{s,u}(k,m)}.
\end{equation*}
The conditional probability  is  well defined  and non-zero only for $k\leq j\leq m\leq N$, and then we have
\begin{equation}\label{cond2}
\Pr(\xi_t=j|\xi_s, \xi_u)=p_{j-\xi_s,\xi_u-\xi_s}(A+\xi_s+t/2,-s+t/2-A-\xi_s, \xi_u+A-t/2+u).
\end{equation}
Indeed, from \eqref{M-P} it follows that
\begin{multline*}
\frac{p_{j-k,N-k}(a,b,c+\delta)p_{m-j,N-j}(a+j-k+\delta,-a-j+k+\delta,c)}{p_{m-k,N-k}(a+\delta,b+\delta,c)}\\=p_{j-k,m-k}(a,b,a+m-k+2\delta).
\end{multline*}
Taking $a=A+t/2+k$, $b=-A-s+t/2-k$, $c=C-u/2$ and $\delta=(u-t)/2$, we get \eqref{cond2}.

We now use \eqref{cond2}  to compute the two-sided conditional moments.
For fixed $\xi_s=k,\xi_u=m$, we use %
\eqref{EX-WilsonD} with $a=A+k+t/2$, $b=-s+t/2-A-k$, $c=m+A-t/2+u$ and $N=m-k$ to compute $\E(Y|\xi_s,\xi_u)$ where $Y=(\xi_t-k)(2a+\xi_t-k)$. %
Since $\Pr(Y_t=y|\xi_s=k,\xi_u=m)=\Pr(Y+(A+k)(A+k+t)=y)$, reverting to $k=\xi_s$ and $m=\xi_u$, we obtain
\begin{multline*}
\E(Y_t|\xi_s,\xi_u)=\left(A+\xi _s\right) \left(A+t+\xi _s\right)\\+\frac{(t-s) \left(\xi
   _u-\xi _s\right) \left(2 A+u+\xi _s+\xi _u\right)}{u-s}
   \\=\frac{u-t}{u-s}\xi_s(2 A+s+ \xi_s)
   +\frac{t-s }{u-s}\xi_u \left(2
   A+u+\xi_u\right)+ A (A+t)  .
  \end{multline*}
Hence a calculation based on \eqref{Y-ini} gives
\begin{equation}\label{LR-Y}
\E(Y_t|Y_s,Y_u) = \frac{(u-t) Y_s+(t-s) Y_u}{u-s}.
\end{equation}

Next, we use \eqref{VX-WilsonD} with $N=\xi_u-\xi_s$, $a=A+\xi_s+t/2$, $b=-s+t/2-A-\xi_s$, and $c=\xi_u+A-t/2+u$ to compute
 the conditional variance:
\begin{multline*}
\Var(Y_t|\xi_s,\xi_u)\\=C_{t,s,u} (\xi_u-\xi_s)(\xi_u-\xi_s+u-s)(2A+s+\xi_s+\xi_u) (2 A+u+\xi_s+\xi_u),
\end{multline*}
where
$$
C_{t,s,u}=\frac{(t-s) (u-t)}{(u-s+1)(u-s)^2}\;.
$$
This gives
\begin{equation}\label{qV-Y}
\Var(Y_t|Y_s,Y_u)=\frac{(u-t)(t-s)}{u-s+1}\left(\frac{(Y_u-Y_s)^2}{(u-s)^2}-\frac{u Y_s-s Y_u}{u-s}\right).
\end{equation}

The following summarizes our findings and incorporates them as an appropriate transformation into a quadratic harness.
\begin{theorem}\label{Cor-parameters}
In Case 1, $(Y_t)_{t\in\TT}$ can be transformed into  a quadratic harness  $(Z_t)_t$ on $\TT'=(0,1)$ with covariance \eqref{EQ:cov} and the conditional variance \eqref{EQ:q-Var} with parameters
\begin{eqnarray}
\eta &=& -\frac{A(B+C+N)+C(B+C)+N(N-B)}{\sqrt{N(A+C)(B+C)(A-B+N)(N-1-B-C)}},\label{eta1}\\
\theta &=& -\frac{A(B+C+N)-B(B+C-N)+2CN}{\sqrt{N(A+C)(B+C)(A-B+N)(N-1-B-C)}}, \label{theta1}\\
\sigma=\tau &=& \frac{1}{N-1-B-C}, \label{sigmatau1}\\
\gamma&=&1+2\sqrt{\sigma\tau}=\frac{B+C-N-1}{B+C-N+1}. \label{gamma1}
\end{eqnarray}

In Case 2, %
$(Y_t)_{t\in\TT}$ can be transformed into a quadratic harness  $(Z_t)_t$ on  $\TT'=(0,\infty)$
 with parameters
\begin{eqnarray*}
\eta &=&-\frac{A(B+C+N)+C(B+C)+N(N-B)}{\sqrt{N(A+C)(B+C)(A-B+N)(B+C-N+1)}},\\
\theta &=& \frac{A(B+C+N)-B(B+C-N)+2CN}{\sqrt{N(A+C)(B+C)(A-B+N)(B+C-N+1)}},\\
\sigma=\tau &=&\frac{1}{B+C-N+1},
\end{eqnarray*}
and $\gamma=1-2\sqrt{\sigma\tau}=(B+C-N-1)/(B+C-N+1)$.
\end{theorem}
\begin{proof}
We shall use Proposition \ref{P-Mobius}. Since Case 1 and Case 2 differ in some details, in order to treat them in a unified way we adopt   the convention that $\eps=1$ refers to Case 1, and $\eps=-1$ refers to  Case 2. We set
\newcommand{\MMM}{M}
\begin{equation}\label{M}
\MMM=\eps(N-B-C)^{-1}\sqrt{\frac{N(A+C)(B+C)(A-B+N)}{\eps(N-1-B-C)}},
\end{equation}
noting that the expression under the radical is positive in both cases.
Let $\psi=A+B$, $\delta=\eps(A-C+N)$, so that $\delta-\eps \psi=\eps(N-B-C)>0$
in both cases.
We take
\[
\alpha=\frac{A(B+C)(A+N)+NBC}{B+C-N},\; \beta=\frac{A(B+C)+NC}{B+C-N},
\]
$\theta_0=0, \eta_0=-1$, and define
\begin{equation}\label{Y2X2Z}
X_t=Y_t-\wo Y_t,\ Z_t=m(t) X_{\ell(t)/m(t)},
\end{equation}
with
\begin{equation}
\ell(t)=\frac{t(A-C+N)-\eps(A+B)}{\MMM (N-B-C)} \mbox{ and } m(t)=\frac{\eps- t}{\MMM (N-B-C)} \label{ell-m}
\end{equation}
 as defined in Proposition \ref{P-Mobius}. Then by Proposition \ref{P-Mobius}, $(Z_t)$   is a quadratic harness, and the formulas follow by  calculation.
 First, $$\chi=\beta^2-\alpha=\frac{(A+C) (B+C) N (A-B+N)}{(B+C-N)^2}>0,$$
 so
$$\MMM/\chi=-\eps(B+C-N)/\sqrt{\eps (A+C)(B+C) N (A-B+N)(N-1-B-C)},$$
and
$\sigma=\tau=\eps/(N-1-B-C)$. Then,
$$\eta= \frac{\MMM (2\beta\eps-\delta)}{\chi}=\frac{\eps \MMM}{\chi}\times \frac{N^2+(A-B) N+(A+C) (B+C)}{B+C-N}$$
and
$$\theta=\frac{\MMM (2\beta-\psi)}{\chi}=\frac{\MMM}{\chi}\times\frac{(A-B) (B+C)+(A+B+2 C) N}{B+C-N}.$$

In view of Theorem \ref{T2}, process $(Z_t)$  in both cases can be defined on $(0,\infty)$ so the one-sided conditional moments are automatically of the correct form.
However, we will still need some of the identities, so we give an argument, which we separate into a lemma.
\end{proof}

\begin{lemma}
$(Z_t)$ satisfies (\ref{LR+}-\ref{RV-}) with parameters, $\eta,\theta,\sigma,\tau$ as given in Theorem \ref{Cor-parameters}.
\end{lemma}

\begin{proof}
From \eqref{cond1} we get
\begin{equation*}
\wo({X}_t|{X}_s)=\frac{A-C+N+t}{A-C+N+s}{X}_s,
\end{equation*}
therefore (note that in both cases $t\mapsto \ell(t)/m(t)$ is increasing)
\begin{equation*}
\wo({X}_{\ell(t)/m(t)}|{X}_{\ell(s)/m(s)})=\frac{1-\eps s}{1-\eps t} {X}_{\ell(s)/m(s)}.
\end{equation*}
Hence $\wo(Z_t|Z_s)=Z_s$.

Next, from  \eqref{vcond1} we get
\begin{multline*}
\Var({X}_t|{X}_s)=\Var(Y_t|Y_s)=
\frac{(A-C+N+t)(t-s)}{(A-C+N+s)^2(A-C+N+s-1)}\\
\times\left[{X}_s+C(s-C)+\wo Y_s\right]\left[{X}_s-(A+N)(s+A+N)+\wo Y_s\right],
\end{multline*}
and
\begin{multline*}
\Var({X}_{\ell(t)/m(t)}|X_{\ell(s)/m(s)})=\eps\frac{(1 -\eps  s)^2(s-t)}{(1-\eps t)^2(1+B+C-N-\eps s)}\\
\times\left[{X}_{\ell(s)/m(s)}+\eps\frac{(A+C)(B+C)}{s-\eps}\right]\times\left[{X}_{\ell(s)/m(s)}+\eps\frac{N(A-B+N)}{s-\eps}\right].
\end{multline*}
This gives
\begin{equation*}
\Var(Z_t|Z_s)=m^2(t)\Var({X}_{\ell(t)/m(t)}|X_{\ell(s)/m(s)})=\frac{t-s}{1+\sigma s}\left(\sigma Z_s+\eta Z_s+1\right).
\end{equation*}

In order to prove \eqref{LR+} and \eqref{RV-}, we define
$$\check{X}_t=Y_{-t},\ t\in \check{\TT}:=-\TT.$$
Then $(\check{X}_t)_t$ is Markov,
\begin{equation*}
\Pr(\check{X}_t=(A-t+k)(A+k))=:\check{p}_{t}(k)
\end{equation*}
with
\begin{equation*}
\check{p}_t(k):=p_{k,N}(A-t/2,B-t/2, C+t/2), \; k=0,\dots,N,
\end{equation*}
 and a standard computation shows that
\[
\Pr\left(\check{X}_t=(A-t+n)(A+n)|\check{X}_s=(A-s+k)(A+k)\right)=:\check{p}_{s,t}(k,n),
\]
where
\begin{equation}\label{2.30}
\check{p}_{s,t}(k,n)=p_{n,k}(A-t/2,B-t/2,A-s+k+t/2),\; n=0,\dots, k,
\end{equation}
(we let $\check{p}_{s,t}(k,n)=0$ for all other values $k,n\in\{0,\dots,N\}$). From Proposition \ref{P:M+V} it follows that
\begin{equation}\label{INV-LR}
\E(\check{X}_t|\check{X}_s)=\check{X}_s\frac{A+B-t}{A+B-s} -\frac{A B (t-s)}{A+B-s},
\end{equation}
and
\begin{multline}\label{INVcondvar}
\Var(\check{X}_t|\check{X}_s)\\=\frac{(A+B-t) (s-t) \left(A (A-s)-\check{X}_s\right) \left(B
   (s-B)+\check{X}_s\right)}{(A+B-s)^2 (A+B-s+1)}.
\end{multline}

For $t,u\in \TT$, $t<u$ we have
\begin{multline*}
\wo( {X}_t|X_u)=\wo(Y_t|Y_u)-\wo Y_t=\wo\left(\check{X}_{-t}|\check{X}_{-u}\right)-\wo Y_t
\\=Y_u
\frac{A+B+t}{A+B+u}-\frac{AB(u-t)}{A+B+u}-\wo Y_t,
\end{multline*}
we obtain that
\[
\wo\left({X}_t|{X}_u\right)={X}_u\frac{A+B+t}{A+B+u}.
\]
Therefore
\[
\wo\left({X}_{\ell(t)/m(t)}|{X}_{\ell(u)/m(u)}\right)=\frac{t(u-\eps)}{u(t-\eps)} {X}_{\ell(u)/m(u)},
\]
so \eqref{LR+} holds true. Similarly, using \eqref{INVcondvar} we get
\begin{multline*}
\Var\left({X}_t|{X}_u\right)=\frac{(A+B+t)(u-t)}{(A+B+u)^2(A+B+u+1)}\\
\times\left[A(A+u)-{X}_u-\wo Y_u\right]\left[B(u+B)- {X}_u-\wo Y_u\right],
\end{multline*}
and
\begin{multline*}
\Var\left({X}_{\ell(t)/m(t)}|{X}_{\ell(u)/m(u)}\right)=\frac{t(u-t)(1-\eps u)^2}{u^2(1-\eps t)^2[1-\eps(1+B+C-N)u]}\\
\times\left[\frac{(A+C)Nu}{u-\eps}+ {X}_{\ell(u)/m(u)}\right]
\left[\eps\frac{(B+C)(A-B+N)u}{1-\eps u}+{X}_{\ell(u)/m(u)}\right].
\end{multline*}
Since $\Var(Z_t|Z_u)=m^2(t)\Var\left( {X}_{\ell(t)/m(t)}| {X}_{\ell(u)/m(u)}\right)$, after a computation, we arrive at
\begin{equation*}
\Var(Z_t|Z_u)=\frac{t(u-t)}{u+\tau}\left(\tau\frac{Z^{2}_u}{u^2}+\theta\frac{Z_u}{u}+1\right).
\end{equation*}

\end{proof}

\section{Extending quadratic harness: conditional representation }\label{Sect:EQHC2}

The next two sections are devoted to the extension of  the quadratic harness  from Case 1 from $(0,1)$ to $(0,\infty)$. We  follow the basic idea suggested by the generalized Waring process  (\cite{burrell1988modelling,burrell1988predictive,zografi2001generalized}), which gives rise to the quadratic harness on $(0,1)$. This quadratic harness can be extended to $(0,\infty)$ by representing the generalized Waring process as a negative binomial process with random parameter, and stitching together two such negative binomial processes that share the randomization, as in  \cite{Maja:2009}. Similarly, we   extend the quadratic harness in Case 1 from  $(0,1)$  to $(0,\infty)$ by representing it as a "Markov process with randomized parameter".  This is assisted here by the heuristic that in Case 1 transition probabilities are positive on  $(-\infty,C-A-N)\cup(-A-B,\infty)$ so there is a natural pair of Markov chains to work with. These two Markov chains can be  put together by requesting that they "match" at infinity, so  the randomization is really based on $\Theta=\lim_{t\to\infty}\xi_t$. (It is clear that once we choose the cadlag trajectories for $\xi_t$, the limit exists almost surely.)

 In this section we analyze two such Markov process, and give the law of the parameter that represents process $(\xi_t)_{t\in\TT}$ from Case 1 as a randomized process. We also give the "dual process" which after randomization would give "the second part" of Case 1 chain, that we did not consider in detail.
 In the next section we stitch together  a pair of such processes.
\subsection{The auxiliary family of  quadratic harnesses}
In this section we construct the family of Markov processes $(\xi_t^{(K)})_{t>-A-B}$ which will give Markov process $(\xi_t)$ from Case 1 once the parameter $K$ is selected at random according to the appropriate law. Heuristically, this process  arises as the limit  $C\to-\infty$ of  the process $(\xi_t)$ 
from Case 1 with $N=K$. But for completeness and for clarity how the remaining parameters enter various formulas we go over the basic analytic identities.

For $K=0,1,\dots$ consider a three-parameter family of finitely supported probability measures $\sum_{j=0}^K \pi_{j,K}(a,b)\delta_j$ on $\{0,1,\dots,K\}$ with probabilities
\begin{multline*}\label{Y-cond}
\pi_{j,K}(a,b)=\frac{(a+1-b)_K}{(a+1)_K}\times \frac{(-1)^j(-K)_j (a)_j (b)_j(1+a/2)_j}{j!(a+1-b)_j(a/2)_j(a+1+K)_j}
\\
=\binom{K}{j}\frac{(a-b+j+1)_{K-j}}{(a+2j+1)_{K-j}}\cdot\frac{(b)_j}{(a+j)_j}.
\end{multline*}
The natural ranges for the parameters are $a>-1$,
$0<b<a+1$, $K\in\NN$, but we also allow $K=0$ with a degenerate law $\delta_0$. The fact that these numbers add up to $1$ can be deduced e.g.
from \cite[identity (9.s)]{Askey:1989} by taking the limit as $c\to\infty$, $e=e(c)\to-\infty$ and $d\to\infty$.
However, it is convenient to observe that
\begin{equation}
  \label{p2pi_lim}
  \pi_{j,K}(a,b)=\lim_{c\to-\infty}p_{j,K}(a/2,b-a/2,c).
\end{equation}
We will rely on this relation for quick proofs of the identities we need.

 We will need moments of the related random variable.
\begin{lemma}\label{L-X-K}
If
$$\Pr(X=j(a+j))=\pi_{j,K}(a,b),$$
then $\E(X)=Kb$ and $\Var(X)=K(K+a-b)b$.
\end{lemma}
\begin{proof} This is  recalculated from the limit as $c\to-\infty$ in Proposition \ref{P:M+V}.

\end{proof}

For each value of $K$, there is a Markov process $\xi_t^{(K)}$ based on these probabilities: the process starts with $\xi_{-A-B}^{(K)}=0$ and has transition probabilities
\begin{equation}\label{Y-CM}
\Pr(\xi_t^{(K)}=j|\xi_s^{(K)}=m)=\pi_{j-m,K-m}(2A+2m+t,t-s).
\end{equation}
(It is straightforward to check that these number are non-negative, and that the univariate laws are $\Pr(\xi_t^{(K)}=j)=\pi_{j,K}(2A+t,t+A+B)$.)
\begin{lemma}
The Chapman-Kolmogorov equations hold.
\end{lemma}
\begin{proof}
The proof is based on the following the algebraic identity:
\begin{equation}
\label{pi-alg}
\frac{\pi_{j,K}(a,b)\pi_{k-j,K-j}(a+\delta+2j,\delta)}{\pi_{k,K}(a+\delta,b+\delta)} = p_{j,k}(a/2,b-a/2,a/2+\delta+k),
\end{equation}
where $p_{j,K}(a,b,c)$ are the previous basic probabilities \eqref{p_k}.
 This identity is recalculated  from \eqref{M-P} using \eqref{p2pi_lim}.

This implies Chapman-Kolmogorov equations in the usual way. We also get the conditional laws under bivariate conditioning: for $s<t<u$,
\begin{multline}\label{CK-K}
\Pr(\xi_t^{(K)}=j|\xi_s^{(K)},\xi_u^{(K)})
\\=p_{j-\xi_s^{(K)},\;\xi_u^{(K)}-\xi_s^{(K)}}(A+\xi_s^{(K)}+t/2,-s+t/2-A-\xi_s^{(K)}, \xi_u^{(K)}+A-t/2+u).
\end{multline}
(This laws are of course the same as  \eqref{cond2}.)
\end{proof}

Next, we define the Markov process of our interest and state the relevant moment formulas.
\begin{proposition}\label{GP}
For $t\in\TT=(-A-B,\infty)$, define $Y_t^{(K)}=(A+t+\xi_t^{(K)})(A+\xi_t^{(K)})$. Then
\begin{enumerate}
\item $(Y_t^{(K)})_{t\in\TT}$ is a Markov process.
\item For $t>-A-B$,
\begin{equation}\label{EY-K}
\E(Y_t^{(K)})=A^2+(A+K)t+K(A+B).
\end{equation}
\item
For $-A-B<s<t$,
\begin{equation}\label{cov-K}
\Cov(Y_s^{(K)}, Y_t^{(K)})= K(K+A-B)(s+A+B).
\end{equation}
\item For $-A-B<s<t$,
\begin{equation}\label{LR-K}
\E(Y_t^{(K)}|Y_s^{(K)})=Y_s^{(K)}+(A+K)(t-s),
\end{equation}
\begin{equation}
\Var(Y_t^{(K)}|Y_s^{(K)})=\left[(A+K)(s+A+K) -Y_s^{(K)}\right](t-s).
\end{equation}
\item For $-A-B<s<t<u$, the two-sided conditional moments are
\begin{equation*}\label{LRk}
\E(Y_t^{(K)}|Y_s^{(K)}, Y_u^{(K)})=
 \frac{(u-t) Y_s^{(K)}+(t-s) Y_u^{(K)}}{u-s},
\end{equation*}
\begin{equation*}
\Var(Y_t^{(K)}|Y_s^{(K)}, Y_u^{(K)})=
\frac{(u-t)(t-s)}{u-s+1}\left(\left(\frac{Y_u^{(K)}-Y_s^{(K)}}{u-s}\right)^2-\frac{u Y_s^{(K)}-s Y_u^{(K)}}{u-s}\right).
\end{equation*}
\item for $-A-B<t<u$, the reverse conditional moments are:
\begin{equation*}\label{LR+K}
\E(Y_t^{(K)}|Y_u^{(K)})=Y_u^{(K)}\frac{A+B+t}{A+B+u} -\frac{A B (u-t)}{A+B+u},
\end{equation*}
\begin{equation*}
\Var(Y_t^{(K)}|Y_u^{(K)})=\frac{(A+B+t) (u-t) \left(Y_u^{(K)}-A (A+u)\right) \left(Y_u^{(K)}-B
   (u+B)\right)}{(A+B+u)^2 (A+B+u+1)}.
\end{equation*}
\end{enumerate}
\end{proposition}
\begin{proof}[Proof  of (i)] %
$(Y_t^{(K)})_{t\in\TT}$ is a one-to-one function of $(\xi_t^{(K)})_{t\in\TT}$.
\end{proof}
\begin{proof}[Proof of (ii)]

From Lemma \ref{L-X-K} with $a=2A+t$,  $b=A+B+t$,  writing  $Y_t^{(K)}=X+A(A+t)$ we
get \eqref{EY-K} and
\begin{equation}\label{Vk}
\Var(Y_t^{(K)})=K(K+A-B)(t+A+B).
\end{equation}
(The latter will be needed for the proof of \eqref{cov-K}.)
 \arxiv{Alternatively, we can take the limit $C\to-\infty$ in Lemma \ref{L.2.3}
 with $N$ exchanged to $K$.}
\end{proof}

\begin{proof}[Proof of (iv)] Comparing \eqref{xiK-tr} and \eqref{Y-CM}, in view of \eqref{p2pi_lim}, the  conditional law of  $\xi_t|\xi_s$ converges as $C\to\-\infty$ to the conditional law of  $\xi_t^{(K)}|\xi_s^{(K)}$. Since the formulas for the Markov processes match, the conditional law
$Y_t^{(K)}|Y_s^{(K)}$ is the limit as $C\to -\infty$ of the conditional laws of the process $(Y_t)$ from Case 1 of Section \ref{Sect:QHFNV}.
So we just pass to  the limit in \eqref{cond1} and \eqref{vcond1}.
\end{proof}
\begin{proof}[Proof of (iii)] This formula follows from \eqref{LR-K} and \eqref{Vk}.
\end{proof}
\begin{proof}[Proof of (v)] Since the conditional laws \eqref{CK-K} are the same as \eqref{cond2}, we use \eqref{LR-Y} and \eqref{qV-Y}.
\end{proof}
\begin{proof}[Proof of (vi)]
For $t<u$, and $j\leq n\leq K\leq N$, the reverse conditional laws  are the same:
\begin{equation}\label{RC-laws}
\Pr(\xi_t^{(K)}=j|\xi_u^{(K)}=n)=\Pr(\xi_t=j|\xi_u=n).
\end{equation}
This follows from the fact that two-sided conditional laws and starting points at $t=-A-B$ are the same.

\end{proof}

Of course, $Y_t^{(K)}\to - AB$ as $t\to-A-B$.
It may be more interesting to remark that once we choose a separable version of the process, we have
 $ \tfrac{Y_t^{(K)}}{t}\to A+K$
almost surely and in mean square as $t\to\infty$.
In particular,
\begin{equation}  \label{YK-conv}
\tfrac{Y_t^{(\Theta)}}{t}\to A+\Theta
\end{equation} almost surely and in mean square for any random $\Theta\in\{0,1,\dots\}$.

\begin{proposition}\label{Prop44} %
If $K\geq 1$ then Markov process
$(Y_t^{(K)})_{t\in\TT}$ can be transformed into a quadratic harness $(Z_t)_t$ %
 on $(0,\infty)$ with parameters
\begin{equation*}
\eta= \frac{1}{\sqrt{K(A-B+K)}}, \; \theta=\frac{A-B+2K}{\sqrt{K(A-B+K)}},
\end{equation*}
$\; \sigma=0$, $ \tau=1$, and  $ \gamma=1$.
\end{proposition}
\begin{proof} The simplest way to get this answer  is to use  Case 1 of Theorem \ref{Cor-parameters} with $N=K$,    taking  the limit as $C\to-\infty$ of the quadratic harness   $(-\sqrt{-C} Z_{-t/C})_{t\in(0,-C)}$.

\arxiv{

Alternatively,  use Proposition \ref{GP} and Proposition \ref{P-Mobius}, keeping in mind that transformation $Z_t\mapsto a Z_{t/a^2}$ maps a quadratic harness with parameters
$\eta,\theta,\sigma,\tau,\gamma$ into a quadratic harness with parameters  $\eta/a,\theta a,\sigma/a^2,\tau a^2,\gamma$.
}
\end{proof}

\subsubsection{Conditional representation}
In this section we confirm that  process $(\xi_t)$ from Case 1 of Section \ref{Sect:QHFNV} can be represented as processes $(\xi_t^{(K)})$ with random $K$.

Denote
\begin{equation}\label{Pi4Theta}
\Pi_k(a,c;N)= \frac{(c)_N(a)_k (-N)_k}{(c-a)_Nk!(c)_k}.
\end{equation}
These numbers are probabilities if $c>0$, $a<1-N$, $k=0,\ldots,N$,
 see \cite[(1.s)]{Askey:1989}.
For $k\ge1$,
$$k\Pi_k(a,c;N)=\frac{a N}{a-c-N+1}\Pi_{k-1}(a+1,c+1;N-1),$$
so  if  $\Pr(U=k)=\Pi_k(a,c;N)$ then
\begin{equation}\label{U-mean-var}
\E(U)= \frac{a N}{a-c-N+1} , \; \Var(U)= \frac{a (a-c+1)  (c+N-1)N}{(c+N-a-2) (c+N-a-1)^2}.
\end{equation}

Consider an auxiliary random variable $\Theta$ with values in $\{0,1,\dots,N\}$ such that
\begin{equation}\label{K-cond}
\Pr(\Theta=k)=\Pi_k(A+C,A-B+1;N),\; k=0, ..., N
\end{equation}
 (This law was calculated from  $\Theta=\lim_{t\to\infty}\xi_t$ in Case 1.)

\arxiv{%
\begin{remark} Recall the constraints introduced at the beginning of Section \ref{Sect:Markov}.
In Case 1 with
$A>-1/2$, $B\in(-A,A+1)$, $C<-A-N+1$, the right hand side of \eqref{K-cond} is indeed positive:  $1-B-C>N-1\ge0$,
$A-B+1>0$ and $A+C+N<1$, so  for $k\in\{0,\dots,N\}$, we have %
$$(A+C)_k(-N)_k=-N(A+C)(-N+1)(A+C+1)\dots(-N+k-1)(A+C+k-1)> 0.$$

In Case 2 with $C>A+N$,  the right hand side of \eqref{K-cond}  is negative when $k+N$ is odd.
\end{remark}
}

\begin{proposition}\label{P-cond-rep} If $\Theta$ is random with law  \eqref{K-cond}, and conditionally on $\Theta=K$,  process $(\xi_t^{(\Theta)})$
is a Markov chain with transitions \eqref{Y-CM},  then the unconditional joint laws are the Case 1 laws: for $j_1\leq \ldots\leq j_n$,
\begin{equation}
  \label{Distr}
  \sum_{k=j_n}^N \Pr(\Theta=k)\Pr\left(\xi_{t_1}^{(k)}=j_1,\dots,\xi_{t_n}^{(k)}=j_n\right)
  =\Pr(\xi_{t_1}=j_1,\dots,\xi_{t_n}=j_n),
\end{equation}
where $(\xi_t)$ is the Markov process from Case 1 with parameters $N,A,B,C$.
\end{proposition}
\begin{proof}
Let $\zeta_t=\xi_t^{(\Theta)}$. Then $(\zeta_t)_{t\in\TT}$ is a Markov chain regardless of the law of the randomization $\Theta$. This follows from the fact that in reverse time the transition probabilities $\xi_t^{(K)}|\xi_u^{(K)}$ do not depend on $K$.
\arxiv{Indeed, for $t<u$,
\begin{multline*}
  \Pr(\zeta_t=j|\zeta_u=k,\Theta)=p_{j,k}(A+t/2,B+t/2,A+u-s/2+k)\indyk{\Theta\geq k}\\
=\indyk{\Theta\geq k}  \Pr(\xi_t=j|\xi_u=k).
\end{multline*}
Since $(\zeta_t)$ is a Markov chain conditionally on $\Theta$, and $\Pr(\zeta_t\leq \Theta)=1$,
 for $j_1\leq j_2\leq\dots\leq j_n$ we have
\begin{multline*}
\Pr(\zeta_{t_1}=j_1,\dots,\zeta_{t_n}=j_n)
\\=
\E\left(\indyk{\Theta\geq j_n} \Pr(\zeta_{t_n}=j_n|\Theta)\prod_{r=1}^{n-1}\Pr(\zeta_{t_r}=j_r| \zeta_{t_{r+1}}=j_{r+1},\Theta)\right)
\\
=  \Pr(\zeta_{t_n}=j_n)\prod_{r=1}^{n-1}\Pr(\xi_{t_r}=j_r| \xi_{t_{r+1}}=j_{r+1}).
\end{multline*}
}
To see that joint laws match, we observe that  Markov processes have the same limiting distribution \eqref{K-cond} and the same
reverse transition probabilities, compare  \eqref{RC-laws}.

\arxiv{\subsubsection*{Direct verification that  $\Pr(\zeta_t=j)=\Pr(\xi_t=j)$} We need to verify that
\begin{multline*}
\sum_{k=j}^N \frac{(A-B+1)_N (A+C)_k (-N)_k}{k! (1-B-C)_N(A-B+1)_k}\pi_{j,k}(2A+t,t+A+B)
\\ = p_{j,N}(A+t/2,B+t/2, C-t/2).
\end{multline*}
 Equivalently,
 $$\sum_{k=j}^N \frac{ (a+c)_k (-N)_k}{k! (a-b+1)_k}\pi_{j,k}(2a,a+b)
 = p_{j,N}(a,b, c)\frac{(1-b-c)_N}{(a-b+1)_N}.$$

This boils down to the following identity:
\begin{multline}\label{Sum-Identity}
\sum_{k=j}^N \frac{ (a+c)_k (-N)_k (-k)_j }{k! (2a+1)_k  (2a+1+k)_j}\\ =(-1)^j \frac{  (a-c+1)_N}{(2 a+1)_N }\times \frac{  (a+c)_j (-N)_j}{ (a-c+1)_j (2 a+N+1)_j}.
\end{multline}
}

\arxiv{%
Simplified form,  see \eqref{jeden}, \eqref{dwa}, \eqref{trzy} is  $$
\sum_{k=j}^N \frac{ (a+c)_k (-N)_k (-k)_j }{k! (2a+1)_{k+j}} =
(-1)^j \frac{  (a-c+1)_N  (a+c)_j (-N)_j }{(a-c+1)_j(2 a+1)_{N+j} }.
$$
Thus
$$
\sum_{k=j}^N \frac{ (a+c)_k  (-1)^{k+j}}{(N-k)!(k-j)! (2a+1)_{k+j}} =
\frac{  (a-c+j+1)_{N-j}  (a+c)_j }{ (2 a+1)_{N+j} (N-j)!}.
$$
Renaming parameters $a+c\mapsto a, a-c\mapsto b$,
$$
\sum_{k=j}^N \frac{ (a)_k  (-1)^{k+j}}{(N-k)!(k-j)! (a+b+1)_{k+j}} =
\frac{  (b+j+1)_{N-j}  (a)_j }{ (a+b+1)_{N+j} (N-j)!}.
$$
Equivalently,
$$
\sum_{k=j}^N \frac{ (a+j)_{k-j}  (-1)^{k+j}}{(N-k)!(k-j)! (a+j+b+1)_{k}} =
\frac{  (b+j+1)_{N-j}  }{ (a+j+b+1)_{N} (N-j)!}.
$$
Renaming $a+j\mapsto a$
}
\arxiv{$$
\sum_{k=j}^N \frac{ (a)_{k-j}  (-1)^{k+j}}{(N-k)!(k-j)! (a+b+1)_{k}} =
\frac{  (b+j+1)_{N-j}  }{ (a+b+1)_{N} (N-j)!}.
$$
  $(a+b+1)_{N}=(a+b+1)_{j} (a+b+j+1)_{N-j}$ and $b+j+1\mapsto b$ give
$$
\sum_{k=j}^N \frac{ (a)_{k-j}  (-1)^{k+j}}{(N-k)!(k-j)! (a+b)_{k-j}} =
\frac{  (b)_{N-j}  }{ (a+b)_{N-j} (N-j)!}.
$$
Changing the index of summation: $k'=k-j$, $N'=N-j$ and dropping the primes, we get
$$
\sum_{k=0}^N \frac{ (a)_{k}  (-1)^{k}}{(N-k)!k! (a+b)_{k}} =
\frac{  (b)_{N} }{ (a+b)_{N} N!}.
$$
Now rename $a+b$ as $c$,
$$\sum_{k=0}^N \big(^N_k\big)\frac{ (a)_{k}  (-1)^{k}}{ (c)_{k}} =
\frac{  (b-c)_{N} }{ (c)_{N} },
$$
and then "undo" the factorials
$$\sum_{k=0}^N \frac{ (a)_{k}  (-N)_{k}}{ k!(c)_{k}} =
\frac{  (b-c)_{N} }{ (c)_{N} }.
$$
This casts \eqref{Sum-Identity} into   \cite[formula (1.s)]{Askey:1989}.
}
\end{proof}

\subsection{The dual process}
For $A,B,C, N$ as in Case 1
and $K=0,\ldots,N$, we now introduce a dual Markov chain
$(\widetilde \xi_t^{(K)})$ 
with state space $\{K,\dots,N\}$ and time $\widetilde \TT=(A+N-C,\infty)$.
This Markov chain starts  at $\widetilde\xi_{A+N-C}^{(K)}=N$ and jumps down according to the transition matrix  $\widetilde P_{s,t}=\left[\widetilde p_{s,t}(i,j)\right]$ with entries
\begin{equation}\label{Dual-cond}
\widetilde p_{s,t}(i,j)=\pi _{j-K,i-K}(2 A+2 K-t,2 A+i+K-s),
\end{equation}
$i=K,K+1,\dots,N,\; j=K,K+1,\dots, i.$
(The remaining entries of this $(N-K+1)\times (N-K+1)$ matrix are zero.)

In particular, the univariate laws of the dual Markov chain are
$$\Pr(\widetilde \xi_t^{(K)}=j)=\pi _{j-K,N-K}(2 A+2 K-t,A+C+K).$$

To confirm that Markov chain $(\widetilde \xi_t^{(K)})_{t\in\widetilde\TT}$ is well defined, we prove the following.  \begin{lemma}
For $s<t$ in $\widetilde \TT$, the entries of transition matrix $\widetilde P_{s,t}$ are non-negative, and the Chapman-Kolmogorov equation holds, i.e. for $s<t<u$, we have
$$\widetilde P_{s,u}= \widetilde P_{s,t}\times\widetilde P_{t,u}.$$
\end{lemma}
\begin{proof}%
Fix $s<t<u$  and $K\leq k\leq j \leq i \leq N$.
We first establish  an identity that will play the role of \eqref{pi-alg} in this argument.
 Taking the limit $b\to\infty$ in \eqref{M-P} we get
    \begin{equation*}
\frac{ \pi_{j,i}(2a+2\delta,c+a+\delta)\pi_{k,j}(2a, 2a+j+2\delta)}{\pi_{k,i}(2a,a+c+\delta)}=p_{j-k,i-k}(a+k+\delta,-a-k+\delta,c).
  \end{equation*}
We use this identity with $\delta=(u-t)/2$, $a=A+K-u/2$, $c=A+K+i+t/2-s$, and then shift the indexes, replacing $i,j,k$ by $i-K,j-K,k-K$ respectively. This gives
 \begin{multline}\label{dual-bicond}
 \frac{\pi _{j-K,i-K}(2 A+2 K-t,2 A+i+K-s)\pi _{k-K,j-K}(2 A+2 K-u,2 A+j+K-t)}{\pi _{k-K,i-K}(2 A+2 K-u,2 A+i+K-s)}\\
= p_{j-k,i-k}(A+k-t/2,u-t/2-A-k,i+A+t/2-s).
\end{multline}
 From \eqref{dual-bicond}  we deduce the Chapman-Kolmogorov equation, and also we determine the two-sided conditional law $\Pr(\widetilde \xi_t^{(K)}=j|\widetilde \xi_s^{(K)}=i,\widetilde \xi_u^{(K)}=k)$. (We omit the verification that the entries are non-negative.)

\arxiv{Here is a direct verification of the non-negativity of the transition probabilities. After a simplification we get
\begin{equation*}
\widetilde p_{s,t}(i,j)=\binom{i-K}{j-K}\cdot\frac{(-(t-s)-(i-j-1))_{i-j}}{(2A-t+2j+1)_{i-j}}\cdot\frac{(2A-s+K+i)_{j-K}}{(2A-t+K+j)_{j-K}}.
\end{equation*}
We are going to use the fact that if $u\in\widetilde{T}$ and $C<-A-N+1$ (as in Case 1), then $2A-u<-2N+1$.

For $i=j$ the first fraction is 1; for $i>j$, since
\begin{eqnarray*}
&&-(t-s)-(i-j-1)\le-(t-s)\le0,\\
&&2A-t+2j+1\le2A-t+j+i<-2N+1+i+j\le0,
\end{eqnarray*}
$\sign(-(t-s)-(i-j-1))_{i-j}=\sign(2A-t+2j+1)_{i-j}$, and the first fraction is non-negative.

Similarly, if $j=K$ then the second fraction is 1; otherwise
\begin{eqnarray*}
&&2A-s+K+i\le2A-s+i+j-1<-2N+i+j\le0,\\
&&2A-t+K+j\le2A-t+2j-1<-2N+2j\le0.
\end{eqnarray*}
Hence $\sign(2A-s+K+i)_{j-K}=\sign(2A-t+K+j)_{j-K}$, and the second fraction is non-negative.

}
\end{proof}

Next, noting again that the lines $\ell_j(t)=(t-A-j) (A+j)$ do not intersect over $\widetilde\TT$, we define the corresponding Markov process
\begin{equation*}\label{du-Y}
\widetilde Y_t^{(K)}=(t-A-\widetilde\xi^{(K)}_t) (A+\widetilde\xi^{(K)}_t)
\end{equation*}

We will need formulas for the absolute moments.
\begin{lemma}
\begin{equation}\label{E-tilde-Y|K}
\E\left( \widetilde Y_t^{(K)}\right)= (A+K) t-(A+K) (A+N)-C (N-K),
\end{equation}
and for $s<t$,
\begin{equation}\label{Var-tilde-Y|K}
\Cov\left(\widetilde Y_s^{(K)},\widetilde Y_t^{(K)}\right)=(K + A + C)(N - K)  (N+A-C-s ).
\end{equation}
\end{lemma}
\begin{proof}
For the mean and variance, we use  Lemma \ref{L-X-K} with $a=2A+2K-t$, $b=A+C+K$ and with $K$ there replaced by $N-K$. Then $\widetilde \xi_t^{(K)}=K+\xi$, where $\xi$ is a random variable such that $X=(a+\xi)\xi$. So
$\widetilde Y_t^{(K)}= -X-(A+K)(A+K-t) $,  and we get both \eqref{E-tilde-Y|K} and the formula for the variance that matches \eqref{Var-tilde-Y|K} when $s=t$.

Next, we apply  Lemma \ref{L-X-K} to  the conditional law \eqref{Dual-cond}. Here
$a=2A+2K-t$, $b=2A+\widetilde \xi_s^{(K)}+K-s$, and the value of $K$ in Lemma \ref{L-X-K} should  now be replaced by $\widetilde \xi_s^{(K)}-K$.   So conditionally on $\widetilde \xi_s^{(K)}$ we can represent $\widetilde \xi_t^{(K)}$ as $K+\xi$, where $\xi$ is a random variable representing $X=(a+\xi)\xi$. Thus conditionally on $\widetilde \xi_s^{(K)}$ we can represent  $\widetilde Y_t^{(K)}$ again as
$-X-(A+K)(A+K-t)$. Since the mean of $X$ is $(\widetilde \xi_s^{(K)}-K)(2A+\widetilde \xi_s^{(K)}+K-s)$,
we get
$$\E(\widetilde Y_t^{(K)}|\widetilde Y_s^{(K)})= (\widetilde \xi_s^{(K)}-K)(s-2A-\widetilde \xi_s^{(K)}-K)-(A+K)(A+K-t)$$
$$
=(s-A-\widetilde \xi_s^{(K)})(A+\widetilde \xi_s^{(K)})+(A+K)(A+K-s)-(A+K)(A+K-t). 
$$
Thus
\[
\E(\widetilde Y_t^{(K)}|\widetilde Y_s^{(K)})=\widetilde Y_s^{(K)}+(A+K)(t-s).
\]
This gives the covariance: from \eqref{G:cov-g} %
we deduce that
$$
\Cov(\widetilde Y_s^{(K)},\widetilde Y_t^{(K)})=\Var(\widetilde Y_s^{(K)}).
$$
\end{proof}

Next, we describe how to get the "second half" of the quadratic harness from Theorem \ref{Cor-parameters}.
\begin{proposition}\label{P-QH-dual}
If $\Theta$ has law \eqref{K-cond} then $(\widetilde Y_t^{(\Theta)})_{t\in\widetilde\TT}$ can be transformed into quadratic harness on $(1,\infty)$ with  parameters
  as in Theorem \ref{Cor-parameters}.

\end{proposition}
\begin{proof}[Sketch of the proof]
We use the fact that $(\widetilde \xi_t^{(\Theta)})_{t\in\widetilde\TT}$ has the same distribution as the time reversal of the original process $(\xi_{-t})_{t\in\widetilde\TT}$ so $\widetilde Y_t^{(\Theta)}= -Y_{-t}$.
With $s<t<u$, this implies that
\begin{equation}\label{EtildeY}
\E(\widetilde Y_t^{(\Theta)})=t\frac{A (B+C)+N C}{B+C-N}-\frac{A(B+C)(A+N)+N B C }{B+C-N},
\end{equation}
 \begin{multline}
\Cov(\widetilde Y_s^{(\Theta)}, \widetilde Y_t^{(\Theta)})\\=-\frac{N(A+C) (B+C) (A-B+N) (t-(A+B)) (s-(A-C+N))}{(B+C-N)^2 (B+C-N+1)}.
\end{multline}
see \eqref{m} and  \eqref{cov}. We also get \eqref{EQ:LR} for $(\widetilde Y_t^{(\Theta)})$
  while \eqref{qV-Y} takes the form
 \begin{equation}\label{VarTilde}
\Var(\widetilde Y_t^{(\Theta)}|\widetilde Y_s^{(\Theta)},\widetilde Y_u^{(\Theta)})=
\frac{(u-t)(t-s)}{u-s+1}\left(\frac{(\widetilde Y_u-\widetilde Y_s)^2}{(u-s)^2}+\frac{u \widetilde Y_s-s \widetilde Y_u}{u-s}\right).
\end{equation}
Let $M$ be given by \eqref{M} with $\eps=1$.
With $\widetilde X_t=\widetilde Y_t^{(\Theta)}-\E(\widetilde Y_t^{(\Theta)})$, taking
\begin{equation}
\ell'(t)=\frac{t(A+N-C)-A-B}{M(N-B-C)} \mbox{ and } m'(t)=\frac{t- 1}{M(N-B-C)}  \label{ell-m1p}
\end{equation}
for $t>1$,
we see that
\begin{equation}\label{Z-dual}
Z_t:= m'(t)\widetilde X_{m'(t)/\ell'(t)}
\end{equation}
defines a Markov process on $(1,\infty)$ such that \eqref{EQ:cov} holds.
A longer calculation verifies \eqref{EQ:q-Var}; this follows  from \eqref{VarTilde}, taking into account \eqref{EtildeY}.
(We remark that  Proposition \ref{P-Mobius}
gives a quadratic harness on $(0,1)$ with  parameters $\eta,\theta$
swapped , i.e. $\eta$ is given by \eqref{theta1} and $\theta$ is given by \eqref{eta1}.
This transformation is based on
$$m(t)=\frac{1-t}{M(N-B-C)} \mbox{ and }\ell(t)=\frac{A+N-C-t(A +B)}{M(N-B-C)}.$$
Then time inversion $tZ_{1/t}$ swaps back the parameters $\eta,\theta$ and maps the process onto $(1,\infty)$. The final transformation is the same as the direct application of  \eqref{ell-m1p}, which is how formula \eqref{Z-dual} was "discovered".)

We omit the verification of one-sided conditional moments which will fall into place anyway since $(Z_t)$ 
extends to a quadratic harness on $(0,\infty)$.
\end{proof}

\section{Extending quadratic harness: stitching two processes together}\label{Sect:stitch}

 Our  goal is to   stitch together a  pair of randomized Markov processes into a single process.
 (The plan of this construction is based on \cite{Bryc-Wesolowski-10}.)
 To do so,
 we chose random variable $\Theta$ with distribution \eqref{K-cond}, and
  a pair of Markov chains $(\zeta_t)_{t\in\TT}$ on $\TT=(-A-B,\infty)$ and $(\zeta_t')$ on $\widetilde\TT=(A+N-C,\infty)$
   such that   $(\zeta_t)$ and $(\zeta'_t)$ are $\Theta$-conditionally independent.
   The law of $(\zeta_t)$ is $(\xi_t^{(\Theta)})$, with state space $\{0,\dots,\Theta\}$ and the law of $(\zeta_t')$ is $(\widetilde \xi_t^{(\Theta)})$ with state space $\{\Theta,\Theta+1,\dots,N-\Theta\}$.

 We then define %
\begin{equation*}\label{Z1}
Z= (A+C)N+\Theta(N-B-C),
\end{equation*}
and
two  (Markov) processes
$Y_t=(A+t+\zeta_t)(A+\zeta_t)$ and $Y_t'=(t-A-\zeta_t')(A+\zeta_t')$.
(Recall that the paths of these processes follow a family of straight lines that do not intersect over $\TT$ and $\widetilde\TT$, so these are indeed Markov processes.)
Let $X_t=Y_t-\E( Y_t)$ and $X_t'=Y_t'-\E (Y_t')$
denote their centered versions.
Processes $(X_t)_{t\in\TT}$, $(X_t')_{t\in\widetilde\TT}$ together with random variable $Z$ will be stitched into a quadratic harness $(Z_t)$ on $(0,\infty)$.

Next we describe the transformations we will use.  Let
 \begin{equation}\label{const}
v= \frac{\sqrt{N (N+A-B)(A+C)(B+C)}}{\sqrt{N-1-B-C}},
\end{equation}
see \eqref{M}.
We then can write \eqref{ell-m} with $\eps=1$
as
\begin{equation}
\ell(t)=(t(A+N-C)-A-B)/v \mbox{ and } m(t)=(1- t)/v \label{ell-m1}
\end{equation} for $0<t<1$.
Similarly, we write \eqref{ell-m1p}
as
\begin{equation}
\ell'(t)=\left(t(A+N-C)-A-B\right)/v \mbox{ and } m'(t)=(t- 1)/v  \label{ell-m2}
\end{equation}
for $t>1$.

The corresponding M\"obius transformations %
are
\begin{equation*}\label{bi-Mobius}
\varphi(t):=\ell(t)/m(t)=\frac{t(A-C+N)-A-B}{1- t} \mbox{ and } \varphi'(t):=\ell'(t)/m'(t)=-\varphi(t).
\end{equation*}
These transformations will be used in the proof. %

The stitched process is then given by
 \begin{equation}\label{Z-abstr}
 Z_t=\begin{cases}
 m(t)X_{\varphi(t) },& 0< t <1,\\
Z/v,& t=1,\\
  m'(t)X'_{\varphi'(t)} ,& t >1,\\
 \end{cases}
 \end{equation}


It is   convenient to observe that   $(Z_t)_{t>0}$ is a Markov process.  Indeed, by Proposition \ref{P-cond-rep}, this follows from Markov property of   $(Y_t)$, and from $\Theta$-conditional independence of $(Y_t)$ and  $(Y_t')$.

The main result  of this section is  the following.
\begin{theorem}\label{T2} For $\Theta$ with law  \eqref{K-cond},
Markov process $(Z_t)_{t>0}$ defined by \eqref{Z-abstr} extends  process $(Z_t)_{t\in(0,1)}$   from Case 1 of Theorem  \ref{Cor-parameters} to
a quadratic harness on $(0,\infty)$ with parameters (\ref{eta1}--\ref{gamma1}).

\end{theorem}

\subsection{Proof of Theorem \ref{T2}}
We need to verify a number of properties from Definition \ref{DEF-harn}. These will be handled after we establish some auxiliary formulas.

\subsubsection{Auxiliary moment calculations}
We first check that $\E(Z)=0$, $\Var(Z)=v^2$ so that $\Var(Z_1)=1$. This is a consequence of the following lemma.
\begin{lemma} \label{L.5.2} For $\Theta$ with law  \eqref{K-cond}, %
$$
\E (\Theta)= \frac{(A+C)N}{B+C-N}  ,\; \Var (\Theta)=M^2,
$$
where $M$ is given by \eqref{M} with $\eps=1$.
\arxiv{ That is,
$$\Var (\Theta)=\frac{N (N+A-B)(A+C)(B+C)}{(N-B-C)^2(N-B-C-1)}.$$
}
\end{lemma}
\begin{proof} See \eqref{U-mean-var}.
\end{proof}

\begin{lemma}\label{L.5.3} For $\Theta$ with law  \eqref{K-cond},
if $0\leq m\leq k\leq n\leq  N$, and $s\in\TT$, $u\in\widetilde \TT$  then
using notation \eqref{Pi4Theta},
$$
\Pr(\Theta=k|\zeta_s, \zeta_u')=\Pi_{k-\zeta_s}(2A+\zeta_s+\zeta_u'-u,2A+2\zeta_s+s +1 ; \zeta_u'-\zeta_s).
$$
\end{lemma}
\begin{proof}
The proof consists of careful isolation of factors that depend only on $k$ in the joint distribution
\begin{multline*}
\Pr(\Theta=k,\zeta_s=m,\zeta_u'=n)
\\=
\Pi_k(A+C,A-B+1;N)\pi_{m,k} (2A+s,s+A+B)\pi_{n-k,N-k}(2A+2k-u,A+C+k)
\\= const_{N,m,n}  \frac{(2A+n+m-u)_{k-m}(-(n-m))_{k-m}}{(k-m)!(2A+2m+s+1)_{k-m}}.
\end{multline*}
(Here  $const_{N,m,n}$ stands for a constant depending on $N,m,n$ and independent of $k$.)
Details are omitted. %
\arxiv{
\begin{proof}Here is a more detailed verification of Lemma \ref{L.5.3}. Isolating the factors that depend on $k$ in $\Pr(\Theta=k,\zeta_s=m,\zeta_u'=n)$, we obtain
\begin{multline*}
\Pr(\Theta=k,\zeta_s=m,\zeta_u'=n)=
\\
\Pi_k(A+C,A-B+1;N)\pi_{m,k} (2A+s,s+A+B)\pi_{n-k,N-k}(2A+2k-u,A+C+k)
\\= const_{N,m,n}\times(-1)^k\binom{N}{k}\binom{k}{m}\binom{N-k}{n-k}\\
\times\frac{(A-B+k+1)_{N-k}(A+C)_k(A-B+m+1)_{k-m}(A+C+k)_{n-k}}{(2A+s+2m+1)_{k-m}(2A-u+k+n)_{n-k}}.
\end{multline*}
(Here and further, $const_{N,m,n}$ is not necessarily the same at each appearance.)

Observe that
\begin{eqnarray*}
(A-B+k+1)_{N-k}(A-B+m+1)_{k-m}&=&(A-B+m+1)_{N-m},\\
(A+C)_k(A+C+k)_{n-k} &=& (A+C)_n,\\
\frac{1}{(2A-u+k+n)_{n-k}} &=& \frac{(2A-u+n+m)_{k-m}}{(2A-u+n+m)_{n-m}}.
\end{eqnarray*}
Since
\begin{equation*}
(-1)^k\binom{N}{k}\binom{k}{m}\binom{N-k}{n-k}=const_{N,m,n}(-1)^k\frac{1}{(k-m)!(n-k)!},
\end{equation*}
and
\[
(-(n-m))_{k-m}=const_{N,m,n}(-1)^k\frac{1}{(n-k)!},
\]
we arrive at the formula
\begin{multline*}
\Pr(\Theta=k,\zeta_s=m,\zeta_u'=n)=\\
const_{N,m,n}\frac{(2A-u+n+m)_{k-m}(-(n-m))_{k-m}}{(k-m)!(2A+s+2m+1)_{k-m}}.
\end{multline*}
Comparing it with \eqref{Pi4Theta}, we get the conclusion of the lemma.
\end{proof}
}

\end{proof}
We will need the first two conditional moments.
\begin{corollary}\label{C-Theta-2cond} For $\Theta$ with law  \eqref{K-cond},
\begin{equation}\label{E(Theta|F)}
\E(\Theta|Y_s, Y_u')= \frac{Y_s+ Y'_u}{u+s}-A,
\end{equation}
 and
 \begin{equation}\label{Var(Theta|Y)}
 \Var(\Theta|Y_s, Y_u')=
\frac{ s^2Y'_u-u^2
   Y_s+su (Y'_u-Y_s) +(Y_s+Y'_u)^2}{(s+u-1) (s+u)^2}.
 \end{equation}

\end{corollary}
\begin{proof}  %
For fixed $\zeta_s=m$, $\zeta'_u=n$,  Lemma \ref{L.5.3} gives
\[
\E(\Theta|\zeta_s=m,\zeta'_u=n)=m+\E(U),
\]
where $\Pr(U=k)=\Pi_k(2A+m+n-u,2A+2m+s+1;n-m)$, so from \eqref{U-mean-var} we get
\begin{multline*}
\E(\Theta|\zeta_s,\zeta_u')=\zeta_s-\frac{(2A+\zeta_s+\zeta_u'-u)(\zeta_u'-\zeta_s)}{s+u}
\\= \frac{1}{s+u}\left(\zeta_s(2A+s+\zeta_s) +\zeta_u' (u-2A-\zeta_u')
\right).
\end{multline*}
This gives \eqref{E(Theta|F)}.

Using \eqref{U-mean-var} again, we get
\begin{multline*}
\Var(\Theta|\zeta_s,\zeta_u')\\= \frac{(\zeta_s - \zeta_u')(\zeta_s + \zeta_u' + 2 A + s)  (\zeta_s + \zeta_u' + 2 A -
   u) (s + u - \zeta_u' + \zeta_s) }{(s + u)^2(s + u - 1)}.
 \end{multline*}
This gives \eqref{Var(Theta|Y)}.

\end{proof}

\subsubsection{Covariance of $(Z_t)$}
\begin{lemma}\label{P2.4}
The stitched process $(Z_t)$ has covariance \eqref{EQ:cov}.
\end{lemma}

\arxiv{
\begin{remark} This should hold true for any law of randomization $\Theta$ when we write the
 conversion  \eqref{Z-abstr} by appropriate transformations that depend on the first two moments of  $\Theta$.
\end{remark}
}

\begin{proof} %
From the transformations  \eqref{Y2X2Z} and \eqref{Z-dual}  exhibited in the proofs of Theorem \ref{Cor-parameters} and Proposition \ref{P-QH-dual}, we see that  the covariance is as required for $0\leq s<u<1$ and   for $1<s<u$, so by time-reversibility argument it remains only to consider the case $s\leq 1<u$.

Since $\lim_{t\to\infty}Y_t/t= A+\Theta $, see \eqref{YK-conv}, we  get  $Z_1=\lim_{s\to 1-}Z_s$  in mean square. %
Therefore, we only need to consider the covariance for  $s<1<u$.
Denote
\begin{equation}\label{rprmt}
s'=\varphi(s)=\ell(s)/m(s) ,\; u'=\varphi'(u)=\ell'(u)/m'(u).
\end{equation}
From  \eqref{Z-abstr}
 we get
$$
\cov(Z_s,Z_u)=m(s)m'(u) \cov(Y_{s'},Y'_{u'})  \,.$$
By conditional independence, from \eqref{G:cov-g} (used with $X=Y_{s'}$, $Y=Y'_{u'}$, $U=\Theta$)
we have
$$\cov(Y_{s'},Y'_{u'})=\cov(\E(Y_{s'}|\Theta),\E(Y'_{u'}|\Theta)).$$
So from  \eqref{EY-K} and \eqref{E-tilde-Y|K}  we get
$$\cov(Y_{s'},Y'_{u'})=(u'-A-N+C)(s'+A+B)  \Var(\Theta).$$
By Lemma \ref{L.5.2}
\begin{multline*}
 \cov(Z_s,Z_u) %
 =M^2 m(s)m'(u) (u'-A-N+C)(s'+A+B)
 \\=M^2 m(s)m'(u)\left[\varphi'(u)-A-N+C\right]\left[\varphi(s)+A+B\right]
 \\
 =M^2 \left[\ell'(u)-(A+N-C)m'(u)\right]\left[\ell(s)+(A+B)m(s)\right].
\end{multline*}
Now we notice that  \eqref{ell-m}, see also \eqref{ell-m1},  %
 gives $\ell(s)+(A+B)m(s)=s/M$, and similarly  \eqref{ell-m1p}, see also \eqref{ell-m2},  %
 gives $\ell'(u)-(A+N-C)m'(u)=1/M$.
 Therefore,
$\cov(Z_s,Z_u)=s$ and  \eqref{EQ:cov} holds.
\end{proof}

\subsubsection{Harness property}

  \begin{lemma}\label{P2.3} Suppose    that  the law of $\Theta$ is \eqref{K-cond}.
Then \eqref{Z-abstr}  defines a harness on $(0,\infty)$.
\end{lemma}
\begin{proof}
The transformations \eqref{Y2X2Z} and \eqref{Z-dual} used in the proofs of Theorem \ref{Cor-parameters} and Proposition \ref{P-QH-dual}, show that \eqref{EQ:LR} holds for $s<t<u<1$ and for $1<s<t<u$.

To end the proof, we only need to verify \eqref{EQ:LR} for $s<t=1<u$.
Indeed, if  we have this case, then the case $0<s<t<1<u$,  is handled from Markov property as
$\E(Z_t|Z_s,Z_u)=\E(\E(Z_t|Z_s,Z_1)|Z_s,Z_u)=\frac{1-t}{1-s}Z_s+\frac{t-s}{1-s} \E(Z_1|Z_s,Z_u)$. The other case $0<s<1<t<u$ is handled similarly (or by time inversion).
Finally, the cases  $1=s<t<u$  and $s<t<u=1$ are the limits of cases $0<s<1<t<u$ and $0<s<t<1<u$, respectively. 

To prove  \eqref{EQ:LR} for  $s<t=1<u$,
we use notation \eqref{rprmt}.
The joint distribution  $Z_s,Z_1,Z_u$ is determined from the joint distribution of $Y_{s'},\Theta,Y'_{u'}$.
To verify harness property, we notice  that Corollary \ref{C-Theta-2cond} implies that $\E\left(Z_1|Y_{s'},Y'_{u'}\right) $ is a linear function of $Y_{s'},Y'_{u'}$, so it is also a linear function of $Z_s,Z_u$. Since by Lemma \ref{P2.4} the covariance of $(Z_t)$ is \eqref{EQ:cov}, this determines the coefficients of the linear regression, and \eqref{EQ:LR} follows.

\end{proof}

\subsubsection{Conditional variance}

\begin{lemma}\label{L6.1}
If $\Theta$ has law \eqref{K-cond}, then \eqref{EQ:q-Var} holds for $t=1$.
\end{lemma}
\begin{proof}
Fix $0<s<1<u$. Using notation \eqref{rprmt},
 we see that $\Var(Z_1|Z_s,Z_u)=\Var(Z_1|Y_{s'},Y'_{u'})$ is a constant multiple of the right hand side of  \eqref{Var(Theta|Y)} (with $s,u$ exchanged to $s',u'$).
  We do not have to pay attention to
  the deterministic multiplicative constant, say $const_{s,u}$, 
  which is determined uniquely from the covariance of $(Z_t)$.
  So we write
  \begin{multline}\label{VarZ1ZsZu}
  \Var(Z_1|Z_s,Z_u)=const_{s,u}\Big(
   {s'}^2Y'_{u'}-{u'}^2
   Y_{s'}+s'u' (Y'_{u'}-Y_{s'}) +(Y_{s'}+Y'_{u'})^2
  \Big).
  \end{multline}
Next, we use the inverse of the transformation \eqref{Z-abstr}, see \eqref{m}, \eqref{EtildeY} and \eqref{rprmt},
  \begin{multline*}
  Y_{s'}= \frac{Z_s}{m(s)}+\frac{A(B+C)(A+N)+N B C }{B+C-N}+\frac{\ell(s)}{m(s)}\frac{A (B+C)+N C}{B+C-N}
  \\=\frac{v Z_s -(A B+C (A+N) s)}{1-s},
  \end{multline*}
  \begin{multline*}
  Y'_{u'}= \frac{Z_u}{m'(u)}+\frac{\ell'(u)}{m'(u)}\frac{A (B+C)+N C}{B+C-N}-\frac{A(B+C)(A+N)+N B C }{B+C-N}
  \\=\frac{v Z_u -(A B+C (A+N) u)}{u-1}.
  \end{multline*}
 Using these expressions, we re-write the right hand side of   \eqref{VarZ1ZsZu}
  as a deterministic multiple of %
\begin{multline*}
 1+\frac{v^2 }{(A+C) (B+C) N (A-B+N)}\left(\frac{(Z_u-Z_s)^2}{(u-s)^2}+\frac{(uZ_s-sZ_u)^2}{(u-s)^2}\right)
 \\-\frac{((A-B) (B+C)+(A+B+2 C) N)
   v}{(A+C) (B+C) N (A-B+N)}\frac{Z_u-Z_s}{u-s}
   \\-\frac{\left(N^2+(A-B) N+(A+C)
   (B+C)\right) v }{(A+C) (B+C) N (A-B+N)} \frac{uZ_s-sZ_u}{u-s}
  \\ +\frac{2 v^2 }{(A+C) (B+C) N (A-B+N)}\frac{(Z_u-Z_s)(uZ_s-sZ_u)}{(u-s)^2}.\end{multline*}
From \eqref{const}, we see that up to a deterministic factor this quadratic form matches \eqref{EQ:q-Var} with parameters (\ref{eta1}-\ref{gamma1}).
\end{proof}

\subsection{Conclusion of proof}

\begin{proof}[Proof of Theorem \ref{T2}]

By Lemma \ref{P2.4},  the covariance  of our process is \eqref{EQ:cov}. From Lemma \ref{P2.3}  we see that \eqref{EQ:LR} holds for all $0<s<t<u$.
Furthermore from the transformations  \eqref{Y2X2Z} and \eqref{Z-dual}  exhibited in the proofs of Theorem \ref{Cor-parameters} and Proposition \ref{P-QH-dual}, we see that $(Z_t)$ is a quadratic harness on $(0,1)$ and on $(1,\infty)$ with the same parameters.
Since \eqref{EQ:LR} holds  for all $0<s<t<u$, from Lemma \ref{L6.1}  and Lemma \ref{Generic-lemma} we see that \eqref{EQ:q-Var} also holds for all $0<s<t<u$.
This ends the proof.
(Recall that the one-sided conditional moments do not need to be verified for processes on $(0,\infty)$.)
\end{proof}

\subsection*{Acknowledgement}
The topic of this research was initiated by a conversation with Albert Cohen.
We would like to thank  J. Weso\l owski for information about \cite{Maja:2009} and several related discussions. We also thank J. A. Wilson for a discussion on the role of discrete models.
The research of WB was partially supported by NSF
grant \#DMS-0904720.

\appendix

\section{Conversion to "standard form"}  In this section we recall a procedure that transforms (some) Markov processes with linear regressions and quadratic conditional variances into the  quadratic harnesses.
The following is \cite[Theorem 3.1]{Bryc-Wesolowski-09}
specialized to
$\chi =0$,  $\eta =
   \eta _0$, $\theta = \theta _0$, $\sigma
   = 0$, $\tau =1$, $\rho = 0$, $a= M$, $b=M \psi $, $c= M
   \epsilon $, $d= M \delta $.
\begin{proposition}\label{P-Mobius}
Suppose $(Y_t)$ is a (real-valued) Markov process on an open interval $\TT\subset\RR$ such that
\begin{enumerate}
\item $\E(Y_t)=\alpha+\beta t$ for some real $\alpha,\beta$.
\item For $s<t$ in $\TT$, ${\rm Cov}(Y_s,Y_t)=M^2(\psi+s)(\delta+\eps t)$, where    $M^2(\psi+t)(\delta+\eps t)>0$ on the entire interval $\TT$, and that $\delta-\eps\psi>0$.
\item For $s<t<u$,
\begin{equation*}\label{computed-cond-var}
\V(Y_t|Y_s,Y_u)=F_{t,s,u}
\left( \eta_0\frac{uY_s-sY_u}{u-s}+\theta_0 \frac{Y_u-Y_s}{u-s}+\frac{(Y_u-Y_s)^2}{(u-s)^2}\right),
\end{equation*}
where $F_{t,s,u}$ is non-random and $ \theta_0,\eta_0\in\RR$  are such that $\chi:=\alpha
\eta_0+\beta\theta_0+\beta^2> 0$.
\end{enumerate}
Denote $X_t=Y_t-\E(Y_t)$. Then there are two affine functions
$$\ell(t)=\frac{t \delta -\psi }{M (\delta - \epsilon
   \psi) }\mbox{ and } m(t)=
 \frac{1-t \epsilon }{M (\delta -\epsilon  \psi) 
   },$$ and an open interval $\TT'\subset(0,\infty)$ such that
$Z_t:=m(t)X_{\ell(t)/m(t)}$ defines a process $(Z_t)$ on $\TT'$  such that \eqref{EQ:cov} holds and  \eqref{EQ:q-Var} holds with parameters
\begin{eqnarray*}
\eta&=& M \left(\delta  \eta _0+\epsilon  \left(2
   \beta +\theta _0\right)\right)/\chi,   \\
\theta&=&   M \left(2 \beta +\psi  \eta _0+\theta
   _0\right)/\chi, \\
\sigma&=& M^2\eps^2/\chi,\\
\tau&=&\ M^2/\chi,\\
\gamma&=& 1+2\eps \sqrt{\sigma\tau}. 
\end{eqnarray*}
\end{proposition}
\begin{remark} The time domain $\TT'$ is the image of $\TT$ under the M\"obius transformation $t\mapsto (t+\psi)/(\varepsilon t+\delta)$.
\end{remark}

\section{Extension Lemma}
The following technical lemma is used  in Section \ref{Sect:stitch}.
\begin{lemma}[\cite{Bryc-Wesolowski-10}]\label{Generic-lemma}
Suppose a square-integrable Markov harness $\Z=(Z_t)_{t>0}$ is a quadratic harness on $(0,1)$ and on $(1,\infty)$, with the same parameters $\ceta,\ctheta,\sigma,\tau,\gamma$.
  If  $\Var(Z_1|Z_s,Z_u)$ is given by the formula \eqref{EQ:q-Var} with $t=1$, and  with the same parameters $\ceta,\ctheta,\sigma,\tau,\gamma$,  then $\Z$ is a quadratic harness on $(0,\infty)$.
\end{lemma}

\arxiv{For completeness, here is a proof from \cite{Bryc-Wesolowski-10}.

 Denote
\begin{equation}\label{Z2Delta}
\Delta_{s,t}=(Z_t-Z_s)/(t-s),\quad \widetilde \Delta_{s,t}=(tZ_s-sZ_t)/(t-s).
\end{equation}

By time-inversion, it suffices to consider formula \eqref{EQ:q-Var} in the case $s<t<1<u$.  By Markov property,
$$
\Var(Z_t|Z_s,Z_u)=\E\big(\Var(Z_t|Z_s,Z_1)\big|Z_s,Z_u\big)+\Var(\E(Z_t|Z_s,Z_1)|Z_s,Z_u)\,.
$$
Denote the right hand side of \eqref{EQ:q-Var} by $F_{t,su}K(Z_s,Z_u)$. Since $\E(Z_t|Z_s,Z_1)$ is given by \eqref{EQ:LR},
 \begin{multline*}\label{tst}
 \Var(\E(Z_t|\calF_{s,1})|\calF_{s,u})=\frac{(t-s)^2}{(1-s)^2}\Var(Z_1|Z_s,Z_u)
 \\=\frac{(t-s)^2(u-1)}{(1-s)(u(1+\sigma s)+\tau-\gamma s)} K(Z_s,Z_u)\,.
\end{multline*}
Next, we write
$$
\E\big(\Var(Z_t|\calF_{s,1})\big|\calF_{s,u}\big)=
\frac{(1-t) (t-s)}{s \sigma +\tau +1-s \gamma }\E(K(Z_s,Z_1)|Z_s,Z_u).
 $$
Since the coefficient $F_{t,s,u}$ is determined by integrating both sides of \eqref{EQ:q-Var}, to end the proof, it suffices to show that
$\E(K(Z_s,Z_1)|Z_s,Z_u)$ is a constant multiple of $K(Z_s,Z_u)$, and we do not need to keep track of the constants.
So it remains to show that
\begin{equation}\label{KKK}
\E(K(Z_s,Z_1)|Z_s,Z_u)=C_{s,u} K(Z_s,Z_u)
\end{equation}
 for any $s<1<u$.

 We have
\begin{equation}\label{K2Delta}
K(Z_s,Z_t)=1+\ceta\widetilde \Delta_{s,t}+\ctheta \Delta_{s,t}+\sigma \widetilde \Delta_{s,t}^2+\tau  \Delta_{s,t}^2-(1-\gamma)\widetilde \Delta_{s,t} \Delta_{s,t}.
\end{equation}

It is easy to check that \eqref{EQ:LR} implies %
\begin{equation}\label{DDD}
  \E(\Delta_{s,t}|\calF_{s,u})=\Delta_{s,u},\; \E(\widetilde \Delta_{s,t}|\calF_{s,u})=\widetilde\Delta_{s,u}.
\end{equation}
}\arxiv{
From  \eqref{DDD} we get
 \begin{eqnarray*}
 \Var(\Delta _{s,t}|\calF_{s,u})&=& \E(\Delta^2_{s,t}|\calF_{s,u})-\Delta_{s,u}^2,\\
 \Var(\widetilde \Delta _{s,t}|\calF_{s,u})&= &\E(\Delta_{s,t}\widetilde\Delta_{s,t}|\calF_{s,u})-\Delta_{s,u}\widetilde\Delta_{s,u},\\
\cov(\Delta _{s,t},\widetilde \Delta _{s,t}|\calF_{s,u})&= &\E(\Delta_{s,t}\widetilde\Delta_{s,t}|\calF_{s,u})-\Delta_{s,u}\widetilde\Delta_{s,u}.
\end{eqnarray*}
Since $\Var(\Delta _{s,t}|\calF_{s,u})$, $\Var(\widetilde \Delta _{s,t}|\calF_{s,u})$ and $\cov(\Delta _{s,t},\widetilde \Delta _{s,t}|\calF_{s,u})$ are all proportional to $\Var(Z_t|\calF_{s,u})$, see \eqref{Z2Delta},
 we get
\begin{eqnarray*}
 \E(\Delta^2_{s,1}|\calF_{s,u})&=&\Delta_{s,u}^2+\frac{1}{(1-s)^2}\Var(Z_1|\calF_{s,u}) \,,\\
  \E(\widetilde\Delta^2_{s,1}|\calF_{s,u})&=&\widetilde\Delta_{s,u}^2+\frac{s^2}{(1-s)^2}\Var(X_1|\calF_{s,u})\,, \\
  \E(\Delta_{s,1}\widetilde\Delta_{s,1}|\calF_{s,u})&=&\Delta_{s,u}\widetilde\Delta_{s,u}-\frac{s}{(1-s)^2}\Var(Z_1|\calF_{s,u})\,.
\end{eqnarray*}
By assumption (ii), $\Var(Z_1|\calF_{s,u})$ is proportional to $K(Z_s,Z_u)$. Using \eqref{K2Delta}, from these formulas together with   \eqref{DDD} we get
$$
\E(K(Z_s,Z_1)|Z_s,Z_u)=K(Z_s,Z_u)+ \frac{\tau+\sigma s^2+(1-\gamma)s}{(1-s)^2} K(Z_s,Z_u),
$$
which proves \eqref{KKK}.
}
\bibliographystyle{alpha}
\bibliography{Vita,W-D-2010}

\begin{thebibliography}{BMW08}

\bibitem[Ask89]{Askey:1989}
R.~Askey.
\newblock {Beta integrals and the associated orthogonal polynomials}.
\newblock {\em Number Theory (ed. K. Alladi). Lecture Notes in Mathematics},
  1395:84--121, 1989.

\bibitem[BMW07]{Bryc-Matysiak-Wesolowski-04}
W{\l}odzimierz Bryc, Wojciech Matysiak, and Jacek Weso{\l}owski.
\newblock Quadratic harnesses, $q$-commutations, and orthogonal martingale
  polynomials.
\newblock {\em Trans. Amer. Math. Soc.}, 359:5449--5483, 2007.
\newblock arXiv.org/abs/math.PR/0504194.

\bibitem[BMW08]{Bryc-Matysiak-Wesolowski-04b}
W{\l}odzimierz Bryc, Wojciech Matysiak, and Jacek Weso{\l}owski.
\newblock The bi-{P}oisson process: a quadratic harness.
\newblock {\em Ann. Probab.}, 36:623--646, 2008.
\newblock arXiv.org/abs/math.PR/0510208.

\bibitem[Bur88a]{burrell1988modelling}
Q.L. Burrell.
\newblock {Modelling the Bradford phenomenon}.
\newblock {\em Journal of Documentation}, 44(1):1--18, 1988.

\bibitem[Bur88b]{burrell1988predictive}
Q.L. Burrell.
\newblock Predictive aspects of some bibliometric processes.
\newblock In L.~Egghe and R.~Rousseau, editors, {\em Informetrics 87/88:
  Selected Proceedings of the First International Conference on Bibliometrics
  and Theoretical Aspects of Information Retrieval}. Elsevier, 1988.

\bibitem[BW10]{Bryc-Wesolowski-08}
W{\l}odek Bryc and Jacek Weso{\l}owski.
\newblock {A}skey--{W}ilson polynomials, quadratic harnesses and martingales.
\newblock {\em Ann. Probab.}, 38(3):1221--1262, 2010.

\bibitem[BW11a]{Bryc-Wesolowski-09}
W{\l}odek Bryc and Jacek Weso{\l}owski.
\newblock Bridges of quadratic harnesses.
\newblock (submitted), 2011.
\newblock arXiv.org/abs/0903.0150.

\bibitem[BW11b]{Bryc-Wesolowski-10}
W{\l}odek Bryc and Jacek Weso{\l}owski.
\newblock Stitching pairs of {L}\'evy processes into martingales.
\newblock In preparation, 2011.

\bibitem[Jam09]{Maja:2009}
Maja Jamio{\l}kowska.
\newblock Bi-{P}ascal process -- definition and properties.
\newblock Master's thesis, Warsaw University of Technology, (in Polish) 2009.

\bibitem[MY05]{Mansuy-Yor04}
Roger Mansuy and Marc Yor.
\newblock Harnesses, {L}\'evy bridges and {M}onsieur {J}ourdain.
\newblock {\em Stochastic Process. Appl.}, 115(2):329--338, 2005.

\bibitem[Wes93]{Wesolowski93}
Jacek Weso{\l}owski.
\newblock Stochastic processes with linear conditional expectation and
  quadratic conditional variance.
\newblock {\em Probab. Math. Statist.}, 14:33--44, 1993.

\bibitem[Wil80]{Wilson:1980}
J.A. Wilson.
\newblock Some hypergeometric orthogonal polynomials.
\newblock {\em SIAM Journal on Mathematical Analysis}, 11:690, 1980.

\bibitem[ZX01]{zografi2001generalized}
M.~Zografi and E.~Xekalaki.
\newblock {The generalized Waring process}.
\newblock In {\em Proceedings of the 5th Hellenic-European Conference on
  Computer Mathematics and its Applications, Athens, Greece}, pages 886--893,
  2001.

\end{thebibliography}
\end{document}